\documentclass[reqno,12pt]{amsart}
\theoremstyle{plain}
\newtheorem{thm}{Theorem}[section]
\newtheorem{cor}[thm]{Corollary} 
\newtheorem{lemma}[thm]{Lemma} 
\newtheorem{prop}[thm]{Proposition}
\newtheorem{defi}[thm]{Definition}

\newtheorem{example}[thm]{Example}

\newtheorem{remark}[thm]{Remark}
\newtheorem{conj}[thm]{Conjecture}

\usepackage{amscd,amssymb,comment,euscript,graphics}
\usepackage{epic}
\usepackage{eepic}
\usepackage[all]{xy}

\newcommand\abar{{\overline a}}

\newcommand\Abar{{\overline A}}

\newcommand\bbar{{\overline b}}

\newcommand\Bbar{{\overline B}}

\newcommand\Bt{{\widetilde B}}

\newcommand\cbar{{\overline c}}

\newcommand\Cbar{{\overline C}}

\newcommand\Cc{{\mathcal{C}}}

\newcommand\Cct{{\widetilde\Cc}}

\newcommand\cl{\operatorname{cl}}

\newcommand\Cpx{{\mathbf C}}

\newcommand\dbar{{\overline d}}

\newcommand\Dbar{{\overline D}}

\newcommand\Dc{{\mathcal{D}}}

\newcommand\diag{\text{\rm diag}}

\newcommand\Eb{{\mathbf E}}

\newcommand\ebar{{\overline e}}

\newcommand\Ebar{{\overline E}}

\newcommand\Ec{{\mathcal{E}}}

\newcommand\eps{\epsilon}

\newcommand\eqdef{{\;\overset{\mbox{\scriptsize def}}{=}}}

\newcommand\fbar{{\overline f}}

\newcommand\Fbar{{\overline F}}

\newcommand\Fc{{\mathcal{F}}}

\newcommand\Fct{{\widetilde\Fc}}

\newcommand\Gc{{\mathcal{G}}}

\newcommand\Gct{{\widetilde\Gc}}

\newcommand\Hbar{{\overline H}}

\newcommand\HEu{{\EuScript H}}                   

\newcommand\Ibar{{\overline I}}

\newcommand\id{{\operatorname{id}}}

\newcommand\ImagPart{{\mathrm{Im}\;}}

\newcommand\Ints{{\mathbf Z}}

\newcommand\Jbar{{\overline J}}

\newcommand\lspan{\mathrm{span}\,}

\newcommand\Mh{{\hat M}}

\newcommand\Nats{{\mathbf N}}

\newcommand\Nbar{{\overline N}}

\newcommand\Oc{{\mathcal{O}}}

\newcommand\oneb{{\mathbf1}}

\newcommand\Pc{{\mathcal{P}}}

\newcommand\phit{{\tilde\phi}}

\newcommand\RealPart{{\mathrm{Re}\;}}

\newcommand\Reals{{\mathbf R}}

\newcommand\restrict{{\upharpoonright}}

\newcommand\Sc{{\mathcal{S}}}

\newcommand\Sf{{\mathfrak S}}

\newcommand\simc{{\;\overset{c}{\sim}\;}}

\newcommand\simcpr{{\;\overset{c'}{\sim}\;}}

\newcommand\SOc{{\mathcal{SO}}}

\newcommand\taut{{\tilde\tau}}

\newcommand\Tc{{\mathcal{T}}}

\newcommand\Tcirc{{\mathbf T}}

\newcommand\Tr{{\mathrm{Tr}}}

\newcommand\Uc{{\mathcal{U}}}

\newcommand\Uct{{\widetilde\Uc}}

\newcommand\Vc{{\mathcal{V}}}

\newcommand\Wc{{\mathcal{W}}}

\begin{document}

\pagestyle{myheadings}

\title{Manifold structure of spaces of spherical tight frames}
\author{Ken Dykema and Nate Strawn}

\address{\hskip-\parindent
Department of Mathematics \\
Texas A\&M University \\
College Station TX 77843--3368, USA}
\email{Ken.Dykema@math.tamu.edu}
\email{natestrawn@neo.tamu.edu}

\thanks{
The first author was supported in part by NSF grant DMS--0300336.
The second author was supported in part by a VIGRE grant from the NSF}

\date{28 September, 2003}

\begin{abstract}
We consider the space $\Fc^\Eb_{k,n}$ of all spherical tight frames of $k$ vectors
in the $n$--dimensional Hilbert space $\Eb^n$ ($k>n$), for $\Eb=\Reals$ or $\Eb=\Cpx$,
and its orbit space $\Gc^\Eb_{k,n}=\Fc^\Eb_{k,n}/\Oc^\Eb_n$ under the obvious action of the group
$\Oc^\Eb_n$ of structure preserving transformations of $\Eb^n$.
We show that the quotient map $\Fc^\Eb_{k,n}\to\Gc^\Eb_{k,n}$ is a locally trivial
fiber bundle (also in the more general case of ellipsoidal tight frames) and
that there is a homeomorphism $\Gc^\Eb_{k,n}\to\Gc^\Eb_{k,k-n}$.
We show that $\Gc^\Eb_{k,n}$ and $\Fc^\Eb_{k,n}$ are real manifolds whenever $k$ and $n$
are relatively prime, and we describe them as a disjoint union of finitely many manifolds (of various
dimensions) when when $k$ and $n$ have a common divisor.
We also prove that $\Fc^\Reals_{k,2}$ is connected ($k\ge4$) and $\Fc^\Reals_{n+2,n}$ is connected,
($n\ge2$).
The spaces $\Gc^\Reals_{4,2}$ and $\Gc^\Reals_{5,2}$ are investigated in detail.
The former is found to be a graph and the latter is the orientable surface of genus $25$.
\end{abstract}

\maketitle

\markboth{Dykema and Strawn}{Manifold structure of frames}

\section{Introduction}

A {\em frame} is a list of vectors $F=(f_i)_{i\in I}$ in a Hilbert space $\HEu$ satisfying
\begin{equation}\label{eq:frame}
A\|v\|^2\le\sum_{i\in I}|\langle v,f_i\rangle|^2\le B\|v\|^2\qquad(v\in\HEu)
\end{equation}
for some constants $0<A\le B$;
the optimal such constants are called the {\em frame bounds} of $F$.
The frame $F$ is finite if the index set $I$ is finite, which implies $\HEu$ is finite dimensional.
An example is an orthonormal basis;
however, in general a frame may have redundancies,
and these are essential in many recent applications of frames (including finite frames) to signal processing
--- see~\cite{GKK}, \cite{CK} and references cited by these.
The frame $F$ is said to be {\em tight} if the constants $A$ and $B$ in~\eqref{eq:frame} can be taken
to be equal to each other.
Some recent references on finite frame theory are \cite{BF}--\cite{HP} and~\cite{SH}.

In this paper, we will consider finite frames, in both the real and complex cases, i.e.\
$\HEu=\Eb^n$ for $\Eb=\Reals$ or $\Eb=\Cpx$.
We will be primarily interested in frames all of whose vectors $f_i$ lie on the unit sphere
of $\Eb^n$, i.e.\ the {\em spherical tight frames}.
(These are also called equal--norm tight frames, uniform tight frames and normalized tight frames
in the literature.)
Our focus will be on the set $\Fc^\Eb_{k,n}$ of all spherical tight frames of $k$ vectors in $\Eb^n$, for $k>n$,
and in particular on the topological questions of connectedness and the manifold structure of $\Fc^\Eb_{k,n}$.

The technical key to our results is to consider the orbit space $\Gc^\Eb_{k,n}=\Fc^\Eb_{k,n}/\Oc^\Eb_n$
for the obvious action of the group of inner--product preserving transformations $\Oc^\Eb_n$ of
the Hilbert space $\Eb^n$.
(Thus, $\Oc_n^\Reals$ is the group of $n\times n$ orthogonal matrices, and $\Oc_n^\Cpx$ is the group
of $n\times n$ unitary matrices.)
It is well known ({\em cf}~\cite{CK} and~\cite{HP})
that $\Gc^\Eb_{k,n}$ can be naturally identified
with the subset of the Grassman manifold of $n$--planes in $\Eb^k$ consisting of projections
all of whose diagonal entries are equal to $n/k$.
We observe that the quotient map $\Fc^\Eb_{k,n}\to\Gc^\Eb_{k,n}$ is a locally trivial
fiber bundle (with fibers $\Oc^\Eb_n$).
In fact, we treat a more general case of ellipsoidal tight frames --- see~\S\ref{sec:eqc}.
An important consequence is that $\Gc^\Eb_{k,n}$ and $\Gc^\Eb_{k,k-n}$ are homeomorphic.

For $n\ge2$, since $\Fc^\Eb_{n+1,1}$ and hence also $\Gc^\Eb_{n+1,1}$ is easy to describe,
the homeomorphism $\Gc^\Eb_{n+1,1}\to\Gc^\Eb_{n+1,n}$ allows us to analyze the space
$\Fc^\Eb_{n+1,n}$ of all spherical tight frames of $n+1$ vectors in $\Eb^n$.
In the real case, we thereby reprove the result~\cite{GKK} that all such frames
are equivalent to each other
if one allows orthogonal transformations of $\Reals^n$ and negating some vectors;
we also prove the analogous result in the complex case.
Finally, we use these techniques to write down explicitly a prototypical example
of a spherical tight frame of $n+1$ vectors in $\Reals^n$, from which all frames
in $\Fc^\Reals_{n+1,n}$ and $\Fc^\Cpx_{n+1,n}$ can be obtained.

Both $\Fc^\Eb_{k,n}$ and $\Gc^\Eb_{k,n}$ are real algebraic sets.
By classical results of Whitney~\cite{W}, each of these
can, therefore, be written as a disjoint union of finitely many manifolds.
We explicitly describe such a decomposition.
When $k$ and $n$ are relatively prime, we show that $\Gc^\Eb_{k,n}$ is itself a real analytic manifold,
and, therefore, so is $\Fc^\Eb_{k,n}$.
When $n$ and $k$ are not relatively prime, $\Gc^\Eb_{k,n}$ is written as a disjoint union of manifolds,
corresponding to block diagonal decompositions of projections.
We get a similar description of $\Fc^\Eb_{k,n}$.
In particular, we say a tight frame $F=(f_1,\ldots,f_k)$ for $\Eb^n$ is {\em orthodecomposable} if
the vectors in $F$ can be partitioned into proper sublists which form tight frames for
orthogonal subspaces of $\Eb^n$.
(See Definition~\ref{def:orthodecomp}.)
Let $\Mh^\Eb_{k,n}$ be the set of spherical tight frames in $\Fc^\Eb_{k,n}$ that are not orthodecomposable.
Then $\Mh^\Eb_{k,n}$ is a nonempty manifold, and $\Fc^\Eb_{k,n}$ is the union of $\Mh^\Eb_{k,n}$
together with other manifolds (of lower dimension) corresponding to orthodecomposability according
to certain partitions.

Another consequence of Whitney's results~\cite{W} is that
$\Fc^\Eb_{k,n}$ and $\Gc^\Eb_{k,n}$ have
only finitely many connected components.
By considering the rearrangement of chains in $\Reals^2$, we prove that the space
$\Fc^\Reals_{k,2}$ of tight spherical frames of $k$ vectors in $\Reals^2$
is connected for all $k\ge4$, and from this result we obtain that the space
$\Fc^\Reals_{n+2,n}$ of real tight spherical frames with two redundant vectors
is connected, for all $n\ge2$.

About half of the length of this paper is occupied with detailed consideration
of two examples:  $\Gc^\Reals_{4,2}$ and $\Gc^\Reals_{5,2}$ (the latter of which is homeomorphic
to $\Gc^\Reals_{5,3}$).
We find that $\Gc^\Reals_{4,2}$ is a graph with twelve vertices and twenty--four edges,
and $\Gc^\Reals_{5,2}$ is the orientable surface of genus $25$.
Similar techniques should permit the description of $\Gc^\Reals_{k,2}$ for larger $k$, though with
considerably more work.
These examples inspired our results on the manifold structure of
$\Gc^\Eb_{k,n}$ for general $k$ and $n$.

The organization of this paper is as follows.
In~\S\ref{sec:eqc}, we show that the quotient map $\Fc^\Eb_{k,n}\to\Gc^\Eb_{k,n}$
is a locally trivial fiber bundle, and the analogous result for ellipsoidal tight frames.
In~\S\ref{sec:1red}, we describe $\Gc^\Eb_{n+1,n}$ and give a concrete example of $F\in\Fc^\Reals_{n,n+1}$.
In~\S\ref{sec:mfld}, we prove $\Gc^\Eb_{k,n}$ and $\Fc^\Eb_{k,n}$
are manifolds when $k$ and $n$ are relatively prime,
and more generally, we write any $\Gc^\Eb_{k,n}$ and $\Fc^\Eb_{k,n}$ as a disjoint union of finitely many
manifolds.
In~\S\ref{sec:G42}, we elucidate $\Gc^\Reals_{4,2}$,
and in~\S\ref{sec:G52} we show $\Gc^\Reals_{5,2}$ is the orientable surface of genus $25$.
In~\S\ref{sec:connected}, we show $\Fc^\Eb_{k,n}$ is connected if and only if $\Fc^\Eb_{k,k-n}$
is connected, and we prove that $\Fc^\Reals_{k,2}$ (for $k\ge4$) and $\Fc^\Reals_{n+2,2}$
(for $n\ge2$) are connected.

\section{Equivalence classes of ellipsoidal tight frames}\label{sec:eqc}

Let $F=(f_1,\ldots,f_k)$ be an ordered frame of $k$ vectors in $\Eb^n$, where $\Eb=\Reals$ or $\Eb=\Cpx$.
Associated to  $F$ is its {\em synthesis operator}
$\Eb^k\to\Eb^n$, defined by
\begin{equation}\label{eq:synth}
\left(\begin{smallmatrix} c_1 \\ \vdots \\ c_k \end{smallmatrix}\right)\mapsto\sum_{j=1}^kc_jf_j.
\end{equation}
The matrix of this operator is thus the $n\times k$ matrix whose columns are the vectors $f_1,\ldots,f_k$ in
this order, and we will identify $F$ with this matrix; thus we also use the notation $F$ for
the synthesis operator~\eqref{eq:synth} itself.
The {\em analysis operator} is the adjoint $F^*:\Eb^n\to\Eb^k$, given by
\[
F^*(v)=\left(\begin{smallmatrix} \langle v,f_1\rangle \\ \vdots \\ \langle v,f_k\rangle
\end{smallmatrix}\right).
\]

Suppose $F$ is a tight frame with frame bound $B$.
Then $B^{-1/2}F^*:\Eb^n\to\Eb^k$ is an isometry.
By a dimensionality argument, there is $U\in\Oc^\Eb_k$ such that
\begin{equation}\label{eq:FWU}
F=B^{1/2}W_{n,k}U,
\end{equation}
where $W_{n,k}=(I_n|0_{n,k-n})$ is the $n\times k$ matrix having $1$
in each $(i,i)$th position and zeros elsewhere.
Conversely, whenever $S:\Eb^n\to\Eb^k$ is an isometry, $F=B^{1/2}S^*$ is a tight frame
of $k$ vectors in $\Eb^n$ having frame bound $B$.

Let $a=(a_1,\ldots,a_n)$ where $a_1\ge a_2\ge\cdots\ge a_n>0$ and consider the ellipsoid
\[
\Ec^\Eb(a)=\{\left(\begin{smallmatrix} v_1 \\ \vdots \\ v_n \end{smallmatrix}\right)\in\Eb^n
\mid\sum_{j=1}^na_j|v_j|^2=1\}.
\]
Letting
\[
D_n(a)=\diag(a_1,\ldots,a_n)\in M_n(\Reals),
\]
we have
\begin{equation}\label{eq:Eca}
\Ec^\Eb(a)=\{v\in\Eb^n\mid\langle D_n(a)v,v\rangle=1\}.
\end{equation}
Let $\Fc^\Eb_k(a)$ denote the set of all ordered tight frames of $k$ vectors that lie on the ellipsoid
$\Ec^\Eb(a)$.
These are the {\em ellipsoidal tight frames} (ETFs) of $k$ vectors on $\Ec^\Eb(a)\subseteq\Eb^n$.
An elementary construction was given in~\cite{DFKLOW} showing that $\Fc^\Eb_k(a)$ is always nonempty.
Let $\Oc^\Eb_n$ act in the usual way on $\Eb^n$ by left multiplication
and let
\[
\Tc^\Eb_n(a)=\{V\in\Oc^\Eb_n\mid V(\Ec^\Eb(a))=\Ec^\Eb(a)\}
\]
be the subgroup of those elements of $\Oc^\Eb_n$ that preserve $\Ec^\Eb(a)$.
From~\eqref{eq:Eca}, we get
\begin{equation}\label{eq:Tca}
\Tc^\Eb_n(a)=\{U\in\Oc^\Eb_n\mid U^*D_n(a)U=D_n(a)\}.
\end{equation}
Then $\Tc^\Eb_n(a)$ acts on $\Fc^\Eb_k(a)$ by left multiplication,
where a frame $F\in\Fc^\Eb_k(a)$ is represented as the $n\times k$ matrix of its synthesis operator,
as described above.
Since the rank of every $F\in\Fc^\Eb_k(a)$ is $n$, this action is free.
We will study the space of orbits of this action.

Let
$\pi:\Fc^\Eb_k(a)\to M_n(\Eb)$
be defined by $\pi(F)=F^*D_n(a)F$.
Since the frame $F\in\Fc^\Eb_k(a)$ consists of vectors lying on the ellipsoid $\Ec^\Eb(a)$, by~\eqref{eq:Eca}
each diagonal entry of $\pi(F)$ is equal to $1$.
Thus $\Tr(\pi(F))=k$.
On the other hand, letting $U\in\Oc^\Eb_k$ be such that $F=B^{1/2}W_{n,k}U$, we have
\begin{equation}\label{eq:piF}
\pi(F)=BU^*W_{n,k}^*D_n(a)W_{n,k}U=BU^*D_k(a)U,
\end{equation}
where
\[
D_k(a)=\diag(a_1,\ldots,a_n,0,\ldots,0)\in M_k(\Reals).
\]
Hence  $\Tr(\pi(F))=B(a_1+\cdots+a_n)$
and the frame bound for $F$ is 
\begin{equation}\label{eq:B}
B=k/(a_1+\cdots+a_n).
\end{equation}
\begin{prop}\label{prop:piT}
Let $F,G\in\Fc^\Eb_k(a)$.
Then $F$ and $G$ lie in the same $\Tc^\Eb_n(a)$--orbit if and only if $\pi(F)=\pi(G)$.
Furthermore, the image of $\pi$ is
\[
\Gc^\Eb_k(a)\eqdef\{R=\tfrac k{a_1+\cdots+a_n}U^*D_k(a)U\mid U\in\Oc^\Eb_k,\,R_{ii}=1,\,(1\le i\le k)\},
\]
where $R_{ii}$ is the $i$th diagonal entry of $R$.
\end{prop}
\begin{proof}
If $F$ and $G$ lie in the same $\Tc^\Eb_n(a)$--orbit, then there is $U\in\Tc^\Eb_n(a)$ such that
$G=UF$.
From~\eqref{eq:Tca}, we get
$\pi(G)=F^*U^*D_n(a)UF=F^*D_n(a)F=\pi(F)$.

On the other hand, suppose $\pi(G)=\pi(F)$.
Let $U,V\in\Oc^\Eb_k$ be such that $F=B^{1/2}W_{n,k}U$ and $G=B^{1/2}W_{n,k}V$.
Using~\eqref{eq:piF} we get $VU^*D_k(a)UV^*=D_k(a)$.
Since $a_n>0$, this yields
\[
VU^*=\left(\begin{smallmatrix}X&0\\0&Y\end{smallmatrix}\right),
\]
where $X\in\Tc^\Eb_n(a)$ and $Y\in\Oc^\Eb_{k-n}$.
Therefore,
\[
G=B^{1/2}W_{n,k}\left(\begin{smallmatrix}X&0\\0&Y\end{smallmatrix}\right)U
=B^{1/2}XW_{n,k}U=XF,
\]
i.e.\ $G$ lies in the same $\Tc^\Eb_n(a)$--orbit as $F$.

The inclusion $\pi(\Fc^\Eb_k(a))\subseteq\Gc^\Eb_k(a)$ is demonstrated in the paragraph
immediately preceding the lemma.
If $R=\frac k{a_1+\cdots+a_n}U^*D_k(a)U\in\Gc^\Eb_k(a)$ for $U\in\Oc^\Eb_k$, then letting $F=B^{1/2}W_{n,k}U$
with $B$ as in~\eqref{eq:B},
we have that $F$ is a tight frame and $R=F^*D_n(a)F$.
If $f_i$ is the $i$th column of $F$, then $1=R_{ii}=\langle D_n(a)f_i,f_i\rangle$,
so $f_i\in\Ec^\Eb(a)$, and therefore $F\in\Fc^\Eb_k(a)$.
\end{proof}

\begin{thm}\label{thm:fiber}
The map $\pi:\Fc^\Eb_k(a)\to\Gc^\Eb_k(a)$ is a locally trivial fiber bundle with fiber $\Tc^\Eb_n(a)$.
\end{thm}
\begin{proof}
From Proposition~\ref{prop:piT} and the freeness of the $\Tc^\Eb_n(a)$--action,
we have that $\pi$ is surjective and, for every $R\in\Gc^\Eb_k(a)$,
$\pi^{-1}(\{R\})$ is homeomorphic to $\Tc^\Eb_n(a)$.
It remains to show local triviality.
For this, it will suffice to find local sections of $\pi$, namely, given $R\in\Gc^\Eb_k(a)$
to find a neighborhood $\Uc$ of $R$ and a continuous map $\mu:\Uc\to\Fc^\Eb_k(a)$ such that $\pi\circ\mu=\id_\Uc$,
because then by Proposition~\ref{prop:piT},
the map $\Tc^\Eb_n(a)\times\Uc\to\Fc^\Eb_k(a)$ given by $(U,S)\mapsto U\mu(S)$
will be a homeomorphism from $\Tc^\Eb_n(a)\times\Uc$ onto $\pi^{-1}(\Uc)$ whose composition with $\pi$
is the projection onto the second component $\Uc$.

Let
\[
\Cc^\Eb_k(a)=\{U^*D_k(a)U\mid U\in\Oc^\Eb_k\}
\]
and let
\[
\sigma:\Oc^\Eb_k\to\Cc^\Eb_k(a)
\]
be $\sigma(U)=U^*D_k(a)U$.
Consider the closed subgroup
\[
\Sc^\Eb_k(a)=\{U\in\Oc^\Eb_k\mid UD_k(a)=D_k(a)U\}
\]
of $\Oc^\Eb_k$ and let $\Sc^\Eb_k(a)\backslash\Oc^\Eb_k$ denote the homogeneous space
of right cosets of $\Sc^\Eb_k(a)$.
The usual quotient map $q:\Oc^\Eb_k\to\Sc^\Eb_k(a)\backslash\Oc^\Eb_k$ is a locally trivial fiber bundle
with fiber $\Sc^\Eb_k(a)$.
The map $\Sc^\Eb_k(a)\backslash\Oc^\Eb_k\to\Cc^\Eb_k(a)$ given by $\Sc^\Eb_k(a)U\mapsto U^*D_k(a)U$
is a homeomorphism.
Hence the map $\sigma$ is a locally trivial fiber bundle with fiber $\Sc^\Eb_k(a)$.

Let
\[
\Cct^\Eb_k(a)=\{S\in\Cc^\Eb_k(a)\mid S_{ii}=\frac{a_1+\cdots+a_n}k,\,(1\le i\le k)\}
\]
and let $\Vc^\Eb_k(a)=\sigma^{-1}(\Cct^\Eb_k(a))$.
The map $r:\Cct^\Eb_k(a)\to\Gc^\Eb_k(a)$ of scalar multiplication by $\frac k{a_1+\cdots+a_n}$
is a surjective homeomorphism.
From the proof of Proposition~\ref{prop:piT},
if $U\in\Vc^\Eb_k(a)$, then letting $F=B^{-1/2}W_{n,k}U$ we have $F\in\Fc^\Eb_k(a)$;
moreover, all elements of $\Fc^\Eb_k(a)$ arise in this way.
Therefore, the map $\rho:\Vc^\Eb_k(a)\to\Fc^\Eb_k(a)$ defined by $\rho(U)=B^{-1/2}W_{n,k}U$
is surjective and continuous, and the diagram
\begin{equation}\label{eq:diag}
\xymatrix{
\Oc^\Eb_k \ar[d]^{\sigma} & \Vc^\Eb_k(a) \ar @{_{(}->} [l] \ar[d]^{\sigma\restrict_{\Vc^\Eb_k(a)}} \ar[r]^{\rho} 
  & \Fc^\Eb_k(a) \ar[d]^{\pi} \\
\Cc^\Eb_k(a) & \Cct^\Eb_k(a) \ar @{_{(}->} [l] \ar[r]^{r} & \Gc^\Eb_k(a)}
\end{equation}
commutes.
Suppose $R\in\Gc^\Eb_k(a)$ and let $R'=\frac{a_1+\cdots+a_n}kR=r^{-1}(R)\in\Cct^\Eb_k(a)$.
There is a neighborhood $\Uc'$ of $R'$ in $\Cct^\Eb_k(a)$ and a local section $\tau:\Uc'\to\Vc^\Eb_k(a)$
of $\sigma$,
which is the restriction to $\Uct\cap\Cct^\Eb_k(a)$
of a local section $\taut:\Uct\to\Oc^\Eb_k$, for some neighborhood $\Uct$ of $R$ in $\Cc^\Eb_k(a)$,
satisfying $\sigma\circ\tau=\id_{\Uc'}$.
Consider the neighborhood $\Uc=r(\Uc')$ of $R$ in $\Gc^\Eb_k(a)$.
Let 
\[
\mu=\rho\circ\tau\circ r^{-1}\restrict_\Uc:\Uc\to\Fc^\Eb_k(a).
\]
Then $\pi\circ\mu=\id_\Uc$.
Hence $\mu$ is a local section of $\pi$.
\end{proof}

\begin{remark}\label{rem:Cinf}\rm
The local section $\taut$ can be taken to be real analytic.
\end{remark}

A consequence of Proposition~\ref{prop:piT}
and Theorem~\ref{thm:fiber} is that $\Gc^\Eb_k(a)$, endowed with the relative topology from $M_k(\Eb)$,
is homeomorphic to the orbit space of the action of $\Tc^\Eb_n(a)$ on $\Fc^\Eb_k(a)$,
endowed with the quotient topology.

\begin{remark}\label{rem:Sk}\rm
If we wish to consider {\em unordered} frames, we should consider the action of
the permutation group $\Sf_k$ on $\Fc^\Eb_k(a)$ by
\[
\Fc^\Eb_k(a)\times\Sf_k\ni(F,\sigma)\mapsto(f_{\sigma(1)},\ldots,f_{\sigma(k)})
=FA_\sigma,
\]
where $A_\sigma$ is the $k\times k$ permutation matrix associated to $\sigma$.
Since this action commutes with the action of $\Tc^\Eb_n(a)$ on $\Fc^\Eb_k(a)$ it descends to the
action
\[
\Gc^\Eb_n(a)\times\Sf_k\ni(R,\sigma)\mapsto A_\sigma^*RA_\sigma
\]
of $\Sf_k$ on $\Gc^\Eb_n(a)$,
and $\Gc^\Eb_n(a)/\Sf_k$ is the orbit space for the action of $\Tc^\Eb(a)$ on the
set of unordered ellipsoidal tight frames.
\end{remark}

\begin{remark}\label{rem:Ek}\rm
Let $\Dc^\Cpx_k=\Tcirc^k$ and $\Dc^\Reals_k=\Dc^\Cpx_k\cap\Reals^k=\{\pm1\}^k$.
If $\zeta=(\zeta_1,\ldots,\zeta_k)\in\Dc^\Eb_k$ and if $F=(f_1,\ldots,f_k)\in\Fc^\Eb_k(a)$,
then setting
\[
F\cdot\zeta=F\diag(\zeta_1,\ldots,\zeta_k)=(\zeta_1f_1,\ldots,\zeta_kf_k),
\]
we have $F\cdot\zeta\in\Fc^\Eb_k(a)$, and this defines an action of the multiplicative
group $\Dc^\Eb_k$ on $\Fc^\Eb_k(a)$.
Since this action commutes with the action of $\Tc^\Eb_n(a)$, it descends to the action of
$\Dc^\Eb_k$ on $\Gc^\Eb_k(a)$ given by
\[
\Gc^\Eb_k(a)\times\Dc^\Eb_k\ni(R,\zeta)\mapsto
\diag(\overline{\zeta_1},\ldots,\overline{\zeta_k})R\,\diag(\zeta_1,\ldots,\zeta_k).
\]
\end{remark}

\medskip

We now specialize to the case of spherical tight frames (STFs), namely when $a_1=\cdots=a_n=1$,
which we will study in the remainder of the paper.
The following corollary restates Proposition~\ref{prop:piT} and Theorem~\ref{thm:fiber}
in this case, and introduces the notation we will use.
As usual, a {\em projection} in $M_k(\Cpx)$ or $M_k(\Reals)$ is a self--adjoint idempotent.
\begin{cor}\label{cor:stf}
Let $\Fc^\Eb_{k,n}$ denote the space of tight frames of $k$ vectors lying on the unit
sphere of $\Eb^n$, and let
\begin{equation}\label{eq:Gkn}
\Gc^\Eb_{k,n}=\{\frac knP\mid P\in M_k(\Eb)\mbox{ a projection of rank }n,\,P_{ii}=\frac nk,\,(1\le i\le k)\}.
\end{equation}
Then the map $\pi=\pi^\Eb_{k,n}:\Fc^\Eb_{k,n}\to\Gc^\Eb_{k,n}$ defined by $\pi(F)=F^*F$ is a surjective,
locally trivial fiber bundle with
fibers $\Oc^\Eb_n$.
Moreover, frames in $\Fc^\Eb_{k,n}$ have the same image under $\pi$ if and only if they
lie in the same orbit of the action of $\Oc^\Eb_n$ on $\Fc^\Eb_{k,n}$.
Hence $\Gc^\Eb_{k,n}$ with the relative topology from $M_k(\Eb)$ is homeomorphic to the space of orbits
of the action of $\Oc^\Eb_n$ on $\Fc^\Eb_{k,n}$, endowed with the quotient topology.
\end{cor}

The following result is now obvious.

\begin{cor}\label{cor:I-P}
If $k,n\in\Nats$ with $k>n$, then there is a homeomorphism $\gamma_{k,n}:\Gc^\Eb_{k,n}\to\Gc^\Eb_{k,k-n}$
given by $\gamma_{k,n}(\frac knP)=\frac k{k-n}(I-P)$.
\end{cor}

\begin{remark}\label{rem:SkEk}\rm
The homeomorphism $\gamma_{k,n}$ intertwines the re--ordering actions of $\Sf_k$ on $\Gc^\Eb_{k,n}$
and $\Gc^\Eb_{k,k-n}$, described in Remark~\ref{rem:Sk}.
Moreover, $\gamma_{k,n}$ intertwines the diagonal actions of $\Dc^\Eb_k$ on $\Gc^\Eb_{k,n}$
and $\Gc^\Eb_{k,k-n}$, described in Remark~\ref{rem:Ek}.
\end{remark}

\section{Frames with one redundant vector}\label{sec:1red}

The spherical tight frames of $n+1$ vectors in $\Reals^n$ are well understood.
Goyal, Kova\v cevi\'c and Kelner proved in~\cite[Thm 2.6]{GKK} that
there is only one of them, up to orthogonal transformations of $\Reals^n$ and the
vector--flipping action of $\Dc^\Reals_{n+1}$ described in Remark~\ref{rem:Ek}.
The homeomorphism $\gamma_{n+1,1}$ of Corollary~\ref{cor:I-P}
yields another proof of this theorem, and of the analogous result for $\Cpx^n$.

\begin{thm}\label{thm:1red}
Let $n\in\Nats$.
Then
\renewcommand{\labelenumi}{(\roman{enumi})}
\begin{enumerate}

\item $\Gc^\Cpx_{n+1,n}$ is homeomorphic to the $n$--torus, $\Tcirc^n$;
morover, the orbit space $\Gc^\Cpx_{n+1,n}/\Dc^\Cpx_{n+1}$ contains only one point;

\item $\Gc^\Reals_{n+1,n}$ has exactly $2^n$ points, $\Gc^\Reals_{n+1,n}/\Sf_{n+1}$ has exactly $[\frac n2]+1$
points, and $\Gc^\Reals_{n+1,n}/\Dc^\Reals_{n+1}$ has only one point.
\end{enumerate}
\end{thm}
\begin{proof}
We prove the complex case~(i), the real case being similar.
The projections of rank~1 in $M_{n+1}(\Cpx)$ having all diagonal entries equal to $1/(n+1)$
are in bijective correspondence with the subspaces of $\Cpx^{n+1}$ of the form $\Cpx v$,
where
\[
v^t=(\tfrac1{\sqrt{n+1}},\tfrac{\zeta_1}{\sqrt{n+1}},\ldots,\tfrac{\zeta_n}{\sqrt{n+1}}),
\]
for $\zeta_1,\ldots,\zeta_n\in\Tcirc$.
This yields the homeomorphism $\Tcirc^n\to\Gc^\Cpx_{n+1,1}$ given by
\[
(\zeta_1,\ldots,\zeta_n)\mapsto(\zeta_{i-1}\overline{\zeta_{j-1}})_{1\le i,j,\le n+1}
\in M_{n+1}(\Cpx),
\]
where we set $\zeta_0=1$.
Since
\[
(\zeta_{i-1}\overline{\zeta_{j-1}})_{1\le i,j,\le n+1}
=\diag(1,\zeta_1,\ldots,\zeta_n)
\left(\begin{matrix}
1&\cdots&1 \\
\vdots&\cdots&\vdots \\
1&\cdots&1
\end{matrix}\right)
\diag(1,\overline{\zeta_1},\ldots,\overline{\zeta_n}),
\]
it is clear that the orbit space $\Gc^\Cpx_{n+1,1}/\Dc^\Cpx_{n+1}$ consists of only one point.
Now the conclusions in~(i) follow from the homeomorphism $\gamma_{n+1,1}$ and Remark~\ref{rem:SkEk}.
\end{proof}

We now write down explicitly a STF, $F$, of $n+1$ vectors in $\Reals^n$.
Thus, all STFs of $n+1$ vectors in $\Reals^n$ are obtained from this one by possibly negating
some vectors and transforming with an element of $\Oc_n^\Reals$,
and all STFs of $n+1$ vectors in $\Cpx^n$ are obtained from $F$ by multiplying the vectors by unimodular
complex numbers and transforming with an element of $\Oc_n^\Cpx$.
\begin{example}\label{ex:def1}\rm
We begin with the frame $F_1=(1,\ldots,1)\in\Fc^\Reals_{n+1,1}$.
The corresponding element of $\Gc^\Reals_{n+1,1}$ is $F_1^*F_1=(n+1)P$, where $P$
is the projection onto the subspace of $\Reals^{n+1}$ spanned by
$w=(1,\ldots,1)^t$.
Applying $\gamma_{n+1,1}$, we get $\frac{n+1}n(I-P)\in\Gc^\Reals_{n+1,n}$, which corresponds
to a frame
\[
F=\sqrt{\frac{n+1}n}\left(I_n\bigg|\begin{smallmatrix}0\\ \vdots\\ 0\end{smallmatrix}\right)V,
\]
where $V\in\Oc^\Reals_{n+1}$ is such that $I-P=V^*\diag(1,\ldots,1,0)V$.
Thus $V$ has rows $v_1,\ldots,v_{n+1}$, where $v_{n+1}^t=\pm w$
and $v_1^t,\ldots,v_n^t$ can be any orthonormal basis for $w^\perp$.
We choose
\begin{align*}
v_1&=\tfrac1{\sqrt2}(1,-1,0,\ldots,0) \\
v_2&=\tfrac1{\sqrt6}(1,1,-2,0,\ldots,0) \\
\vdots \\
v_j&=\tfrac1{\sqrt{j(j+1)}}(\underset{j}{\underbrace{1,\ldots,1}},-j,0,\ldots,0) \\
\vdots \\
v_n&=\tfrac1{\sqrt{n(n+1)}}(1,\ldots,1,-n).
\end{align*}
This yields the frame $F=(f_1,\ldots,f_{n+1})\in\Fc_{n+1,n}$,
where
\begin{alignat*}{2}
f_1^t&=\sqrt{\tfrac{n+1}n}\Big(\tfrac1{\sqrt2},\tfrac1{\sqrt6},\tfrac1{\sqrt{12}},\tfrac1{\sqrt{20}},
 \ldots,\tfrac1{\sqrt{n(n+1)}}\Big) \\
f_2^t&=\sqrt{\tfrac{n+1}n}\Big(\tfrac{-1}{\sqrt2},\tfrac1{\sqrt6},\tfrac1{\sqrt{12}},\tfrac1{\sqrt{20}},
 \ldots,\tfrac1{\sqrt{n(n+1)}}\Big) \\
\vdots \\
f_p^t&=\sqrt{\tfrac{n+1}n}\Big(\underset{p-2}{\underbrace{0,\ldots,0}},
 \tfrac{-(p-1)}{\sqrt{(p-1)p}},\tfrac1{\sqrt{p(p+1)}},\ldots,\tfrac1{\sqrt{n(n+1)}}\Big), \\
\vdots \\
f_{n+1}^t&=\sqrt{\tfrac{n+1}n}\Big(0,\ldots,0,\tfrac{-n}{\sqrt{n(n+1)}}\Big)&&=(0,\ldots,0,-1).
\end{alignat*}
One easliy verifies that all vectors in this frame have the same angle between them:
\[
\langle f_p,f_q\rangle=-1/n,\quad(p\ne q).
\]
\end{example}

\section{Manifold structure}\label{sec:mfld}

Let $\Eb=\Reals$ or $\Eb=\Cpx$.
Let $\odot$ denote the binary operation of coordinate--wise multiplication on $\Eb^k$.
Thus, if $v=(v_i)_{i=1}^k$ and $w=(w_i)_{i=1}^k$, then $v\odot w=(v_iw_i)_{i=1}^k$.
Taking $\odot$ as multiplication makes $\Eb^k$ into the commutative, unital $\Eb$--algebra
that is often denoted $\ell^\infty_{k,\Eb}$

Let $k,n\in\Nats$, $k>n$.
The Grassman manifold
\begin{equation}\label{eq:Grassman}
G^\Eb_{k,n}=\{V^*DV\mid V\in\Oc^\Eb_k\},\qquad D=\diag(\underset{n}{\underbrace{1,\ldots,1}},0,\ldots,0)
\end{equation}
of $n$--planes in $\Eb^k$ is a real analytic submanifold of $M_k(\Eb)$, and
we have
\[
\frac nk\Gc^\Eb_{k,n}=\{P\in G^\Eb_{k,n}\mid P_{ii}=\frac nk,\,(1\le i\le k)\}.
\]
Let
\[
f:G^\Eb_{k,n}\to K_{k,n}\eqdef\{(d_i)_{i=1}^k\in\Reals^k\mid\sum_{i=1}^kd_i=n\}
\]
be the map $f(P)=(P_{ii})_{i=1}^k$, that extracts the diagonal of the projection.
For $P\in G^\Eb_{k,n}$, we denote as usual the differential of $f$ at $P$ by
\begin{equation}\label{eq:DfP}
Df_P:T_PG^\Eb_{k,n}\to T_{f(P)}K_{k,n}.
\end{equation}

Let $\xi_1,\ldots,\xi_k$ be the standard orthonormal basis for $\Reals^k$, (also for $\Cpx^k$).
Given a subset $A\subseteq\{1,\ldots,k\}$, let $E_A=\sum_{i\in A}\xi_i$ and
let $Q_A:\Eb^k\to\Eb^k$ be the projection onto the subspace 
\[
\lspan\{\xi_i\mid i\in A\}
\]
of $\Eb^k$.
Then 
\begin{equation}\label{eq:Qeta}
Q_A(\eta)=E_A\odot\eta,\qquad(\eta\in\Eb^k).
\end{equation}
\begin{defi}\label{def:sigmaT}\rm
Given $T\in M_k(\Eb)$, let $\sigma_T$ be the set of all minimal nonempty subsets $A$ of $\{1,\ldots,k\}$
that satisfy $Q_AT=TQ_A$.
Note that $\sigma_T$ is a partition of $\{1,\ldots,k\}$.
\end{defi}
\begin{lemma}\label{lem:regpt}
Take a projection $P\in G^\Eb_{k,n}$ in the Grassman manifold.
Let $\Wc\subseteq\Reals^k$ be the range of the differential map~\eqref{eq:DfP},
and let $\Wc^\perp$ denote the orthocomplement of $\Wc$ in $\Reals^k$.
Then a basis for $\Wc^\perp$ is
\begin{equation}\label{eq:Wperpbasis}
\{E_A\mid A\in\sigma_P\}.
\end{equation}
Consequently, $P$
is a regular point of $f$ if and only if $PQ_A\ne Q_AP$ for all
proper, nonempty subsets $A\subseteq\{1,\ldots,k\}$.
\end{lemma}
\begin{proof}
Since $P$ is a regular point of $f$ if and only if $\dim(\Wc)=k-1$,
the last statement of the lemma will follow immediately once~\eqref{eq:Wperpbasis} is shown to be a basis
for $\Wc^\perp$.

Take $\Eb=\Cpx$.
Let $P=V^*DV$ be as in~\eqref{eq:Grassman}.
Let $(e_{\iota j})_{1\le\iota,j\le k}$ be the standard system of matrix units for $M_k(\Cpx)$.
A basis for the tangent space $T_PG^\Cpx_{k,n}$ is the list of $2n(k-n)$ vectors
\[
\big(x(\iota,j)\big)_{1\le\iota\le n<j\le k},\quad\big(y(\iota,j)\big)_{1\le\iota\le n<j\le k},
\]
where
\begin{align*}
x(\iota,j)&=\frac d{dt}\bigg|_{t=0}V^*e^{t(e_{j\iota}-e_{\iota j})}De^{t(e_{\iota j}-e_{j\iota})}V
 =V^*(e_{j\iota}+e_{\iota j})V \\
y(\iota,j)&=\frac d{dt}\bigg|_{t=0}V^*e^{-it(e_{\iota j}+e_{j\iota})}De^{it(e_{\iota j}+e_{j\iota})}V
 =iV^*(e_{\iota j}-e_{j\iota})V.
\end{align*}
The $p$th diagonal entries of these are
\[
x(\iota,j)_{pp}=2\,\RealPart(v_{\iota p}\overline{v_{jp}}),\qquad
y(\iota,j)_{pp}=2\,\ImagPart(v_{\iota p}\overline{v_{jp}}),
\]
where $v_{\iota p}$ is the $(\iota,p)$th entry of $V$.
Let $v_\iota$ denote the $\iota$th row of $V$.
Therefore,
\[
Df_P(x(\iota,j))=2\,\RealPart(v_\iota\odot\overline{v_j}),\qquad
Df_P(y(\iota,j))=2\,\ImagPart(v_\iota\odot\overline{v_j}).
\]
Letting $\Vc=P(\Cpx^k)$, we have
\begin{align*}
\Vc&=\lspan\{v_\iota^t\mid1\le\iota\le n\} \\
\Vc^\perp&=\lspan\{v_j^t\mid n<j \le k\}.
\end{align*}
For $u\in\Reals^k$, we therefore have
\begin{equation}\label{eq:uV}
\begin{aligned}
u\in\Wc^\perp\quad&\Leftrightarrow\quad\langle u,v\odot\overline{v'}\rangle=0,
\quad(v\in\Vc,\,v'\in\Vc^\perp) \\
&\Leftrightarrow\quad\langle v',u\odot v\rangle=0,
\quad(v\in\Vc,\,v'\in\Vc^\perp) \\
&\Leftrightarrow\quad u\odot\Vc\subseteq\Vc.
\end{aligned}
\end{equation}
From~\eqref{eq:uV}, we see that $\Wc^\perp$ is a unital subalgebra of $\ell^\infty_{k,\Reals}$.
It is a standard result, and not difficult to show, that all unital subalgebras 
of $\ell^\infty_{k,\Reals}$ are of the form
\begin{equation}\label{eq:linfsa}
\lspan_\Reals\{E_A\mid A\in\sigma\},
\end{equation}
where $\sigma$ is a partition of $\{1,\ldots,k\}$.
But from~\eqref{eq:uV} and~\eqref{eq:Qeta},
\[
E_A\in\Wc^\perp\quad\Leftrightarrow\quad E_A\odot\Vc\subseteq\Vc
\quad\Leftrightarrow\quad Q_A(\Vc)\subseteq\Vc
\quad\Leftrightarrow\quad Q_AP=PQ_A.
\]
This concludes the proof in the case $\Eb=\Cpx$.

The proof in the case $\Eb=\Reals$ is similar, but easier.
Indeed, a basis for the tangent space of $G^\Reals_{k,n}$ is $(x(\iota,j))_{1\le\iota\le n<j\le k}$
and, with $\Vc=P(\Reals^k)$, we find that the range of $Df_P$ is
\[
\Wc=\lspan\{v\odot v'\mid v\in\Vc,\,v'\in\Vc^\perp\}.
\]
Now the proof proceeds as before, beginning with the chain of implications~\eqref{eq:uV}.
\end{proof}

\begin{thm}\label{thm:mfld}
Let $n,k\in\Nats$, $k>n$, with $n$ and $k$ relatively prime.
Then
\renewcommand{\labelenumi}{(\roman{enumi})}
\begin{enumerate}

\item
$\Gc^\Reals_{k,n}$ is a regular, real analytic submanifold
of $M_k(\Reals)$ of dimension 
\[
\dim(\Gc^\Reals_{k,n})=(k-n-1)(n-1);
\]

\item
$\Fc_{k,n}^\Reals$ is a regular, real analytic submanifold of $(S^{n-1})^k$ of dimension
\[
\dim(\Fc^\Reals_{k,n})=(k-\frac n2-1)(n-1);
\]

\item
$\Gc^\Cpx_{k,n}$ is a regular, real analytic submanifold
of $M_k(\Cpx)$ of dimension
\[
\dim(\Gc^\Cpx_{k,n})=2n(k-n)-k+1;
\]

\item
$\Fc_{k,n}^\Cpx$ is a regular, real analytic submanifold of $(S^{2n-1})^k$ of dimension
\[
\dim(\Fc^\Cpx_{k,n})=2n(k-n)+n^2-k+1;
\]

\end{enumerate}
\end{thm}
\begin{proof}
By the proof of Theorem~\ref{thm:fiber} and Remark~\ref{rem:Cinf},~(ii)
will follow from~(i) and~(iv) will follow from~(iii).

We will show that $c=(\frac nk)_{i=1}^k$ is a regular value of $f:G^\Eb_{k,n}\to K_{k,n}$,
for $\Eb=\Reals$ and $\Eb=\Cpx$.
Since $\Gc^\Eb_{k,n}$ is nonempty (see~\cite{GKK}, \cite{HP} or~\cite{DFKLOW}),
by the regular value
theorem, this will imply~(i) and~(iii).
By Lemma~\ref{lem:regpt}, it will suffice to show that if $P\in f^{-1}(c)$, then
$PQ_A\ne Q_AP$ for all proper, nonempty subsets $A$ of $\{1,\ldots,k\}$.
Suppose, to obtain a contradiction, we have $Q_AP=PQ_A$ for some such subset $A$.
Then $PQ_A$ can be viewed as an $|A|\times|A|$ matrix
and is a projection, all of whose diagonal entries are $\frac nk$.
The rank of $Q_AP$ is thus $\frac nk|A|$.
However, since $1\le |A|\le k-1$ and since $n$ and $k$ are relatively prime,
$\frac nk|A|$ cannot be an integer;
this is a contradiction.
\end{proof}

\begin{lemma}\label{lem:dregpt}
Let $k,n\in\Nats$, $k>n$.
Then there is $S\in\Gc^\Reals_{k,n}$ such that $\frac nkS$ is a regular point of the map $f$.
\end{lemma}
\begin{proof}
Let $d$ be the greatest common divisor of $k$ and $n$.
If $d=1$, then it follows from the proof of Theorem~\ref{thm:mfld}
that for every $S\in\Gc^\Reals_{k,n}$,
$\frac nkS$ is a regular point of $f$.
Suppose $d>1$.
Let $k'=k/d$ and $n'=n/d$.
Let $R'\in\Gc^\Reals_{k',n'}$ and let
$R=\diag(R',\ldots,R')$
be the indicated block diagonal $d\times d$ matrix of $k'\times k'$ matrices.
Then $R\in\Gc^\Reals_{k,n}$.
Let $\xi_1,\ldots,\xi_k$ be the standard orthonormal basis of $\Reals^k$ and for $\ell<k$,
identify $\Reals^\ell$
with the usual subspace of $\Reals^k$, having standard orthonormal basis $\xi_1,\ldots,\xi_\ell$.
Let $U\in\Oc^\Reals_d$ be a real orthogonal matrix satisfying
\[
\langle U\xi_1,\xi_j\rangle\ne0,\qquad(1\le j\le d).
\]
Let $W\in\Oc^\Reals_k$ be any real orthogonal matrix satisfying
\[
W\xi_j=\xi_{1+(j-1)k'},\qquad(1\le j\le d).
\]
Let $V=W\Big(\begin{smallmatrix} U & 0 \\ 0 & I_{k-d} \end{smallmatrix}\Big)W^*\in\Oc^\Reals_k$.
Then
\[
V\xi_p=\begin{cases}
\sum_{i=1}^du_{ij}\xi_{1+(i-1)k'},&\quad p=1+(j-1)k',\,1\le j\le d, \\
\xi_p,&\quad p\notin\{1,\,1+k',\,1+2k',\,\ldots,\,1+(d-1)k'\},
\end{cases}
\]
where $u_{ij}$ is the $(i,j)$th entry of $U$.
Let $S=V^*RV$.

In order to show $S\in\Gc^\Reals_{k,n}$, it will suffice to show
\begin{equation}\label{eq:Spp}
\langle S\xi_p,\xi_p\rangle=1,\qquad(1\le p\le k).
\end{equation}
Suppose $p\notin\{1,\,1+k',\,1+2k',\,\ldots,\,1+(d-1)k'\}$.
Then
\[
\langle S\xi_p,\xi_p\rangle=\langle RV\xi_p,V\xi_p\rangle=\langle R\xi_p,\xi_p\rangle=1.
\]
Suppose $p=1+(j-1)k'$.
Then
\begin{align*}
\langle S\xi_p,\xi_p\rangle&=\langle RV\xi_p,V\xi_p\rangle
=\sum_{i=1}^d\sum_{i'=1}^du_{ij}\overline{u_{i'j}}\langle R\xi_{1+(i-1)k'},\xi_{1+(i'-1)k'}\rangle \\
&=\sum_{i=1}^d|u_{ij}|^2\langle R\xi_{1+(i-1)k'},\xi_{1+(i-1)k'}\rangle=\sum_{i=1}^d|u_{ij}|^2=1.
\end{align*}
Thus~\eqref{eq:Spp} is proved.

Consider the relation $\simc$ on $\{1,\ldots,k\}$ defined by
\[
i\simc j\quad\Leftrightarrow\quad\langle S\xi_i,\xi_j\rangle\ne0
\]
and let $\sim$ be the equivalence relation on $\{1,\ldots,k\}$ generated by $\simc$.
We will show that $\sim$ has
only one equivalence class, which by Lemma~\ref{lem:regpt}, is equivalent to
$\frac nkS$ being a regular point of $f$.
Let $\simcpr$ be the relation on $\{1,\ldots,k'\}$ defined by
\[
i\simcpr j\quad\Leftrightarrow\quad\langle R'\xi_i,\xi_j\rangle\ne0.
\]
Since $k'$ and $n'$ are relatively prime,
from the proof of Theorem~\ref{thm:mfld} we have $Q_AR'=R'Q_A$ for all proper, nonempty
subsets $A\subseteq\{1,\ldots,k\}$, and therefore,
we know that
the equivalence relation on $\{1,\ldots,k'\}$
generated by $\simcpr$ has only one equivalence class.
Take $i'\in\{1,\ldots,d\}$ and $s,t\in\{2,\ldots,k'\}$ and set
\[
p=s+(i'-1)k',\qquad q=t+(i'-1)k'.
\]
Then
\begin{align*}
\langle S\xi_p,\xi_q\rangle&=\langle R\xi_p,\xi_q\rangle=\langle R'\xi_s,\xi_t\rangle \\
\langle S\xi_1,\xi_q\rangle&=\sum_{i=1}^du_{i1}\langle R\xi_{1+(i-1)k'},\xi_q\rangle
=u_{i'1}\langle R\xi_{1+(i'-1)k'},\xi_q\rangle=u_{i'1}\langle R'\xi_1,\xi_t\rangle.
\end{align*}
Therefore,
\begin{align*}
s\simcpr t&\quad\implies\quad s+(i'-1)k'\simc t+(i'-1)k', \\
1\simcpr t&\quad\implies\quad 1\simc t+(i'-1)k'.
\end{align*}
We conclude
\begin{equation}\label{eq:1simq}
1\sim q,\qquad(q\in\{1,\ldots,k\}\backslash\{1,\,1+k',\,1+2k',\,\ldots,\,1+(d-1)k'\}).
\end{equation}
Finally, there must be $s\in\{2,\ldots,k'\}$ such that $1\simcpr s$.
Given $i\in\{1,\ldots,d\}$,
let $i'\in\{1,\ldots,d\}$ be such that $u_{i'i}\ne0$ and let
\[
p=1+(i-1)k',\qquad q=s+(i'-1)k'.
\]
Then
\[
\langle S\xi_p,\xi_q\rangle=\langle RV\xi_p,\xi_q\rangle
=\sum_{j'=1}^d u_{j'i}\langle R\xi_{1+(j'-1)k'},\xi_q\rangle
=u_{i'i}\langle R'\xi_1,\xi_s\rangle\ne0.
\]
Thus, we have $p\simc q$.
But from~\eqref{eq:1simq} we have $q\sim 1$, and we conclude $p\sim 1$.
Combined with~\eqref{eq:1simq}, this shows $1\sim r$ for all $r\in\{1,\ldots,k\}$.
\end{proof}

\begin{defi}\label{def:N}\rm
Let $k,n\in\Nats$ with $k>n$ and let $\Eb=\Reals$ or $\Eb=\Cpx$.
Whenever $\sigma$ is a partition of the set $\{1,\ldots,k\}$, let
\[
N^\Eb_{k,n}(\sigma)=\{R\in\Gc^\Eb_{k,n}\mid \sigma_R=\sigma\},
\]
where $\sigma_R$ is as in Definition~\ref{def:sigmaT}.
We will also write simply $N^\Eb_{k,n}$ to denote $N^\Eb_{k,n}(\oneb_k)$,
where $\oneb_k$ is the trivial partition of $\{1,\ldots,k\}$ into one subset.
\end{defi}

\begin{thm}\label{thm:strat}
Let $k,n\in\Nats$ with $k>n$ and let $\Eb=\Reals$ or $\Eb=\Cpx$.
Then
\renewcommand{\labelenumi}{(\roman{enumi})}
\begin{enumerate}

\item
$N^\Eb_{k,n}$ is a nonempty, regular, real analytic submanifold of $M_k(\Eb)$ with
\[
\dim(N^\Eb_{k,n})=\begin{cases}
(k-n-1)(n-1),&\quad\Eb=\Reals, \\
2n(k-n)-k+1,&\quad\Eb=\Cpx.
\end{cases}
\]

\end{enumerate}
Let $d=\gcd(k,n)$ and let $k'=k/d$, $n'=n/d$.
Let $\Pc(k,k')$ be the set of all partitions of the set $\{1,\ldots,k\}$
into subsets whose cardinalities are multiples of $k'$.
Then
\renewcommand{\labelenumi}{(\roman{enumi})}
\begin{enumerate}
\stepcounter{enumi}

\item 
\begin{equation}\label{eq:Gcdecomp}
\Gc^\Eb_{k,n}=\bigcup_{\sigma\in\Pc(k,k')}N^\Eb_{k,n}(\sigma)
\end{equation}
and the sets $(N^\Eb_{k,n}(\sigma))_{\sigma\in\Pc(k,k')}$ are pairwise disjoint;

\item
if $\sigma=\{A_1,\ldots,A_\ell\}\in\Pc(k,k')$ with $|A_i|=m_ik'$,
then $N^\Eb_{k,n}(\sigma)$ is a nonempty, regular, real analytic submanifold of $M_k(\Eb)$
and is real--analytically diffeomorphic to the Cartesian product
\[
\prod_{i=1}^\ell N^\Eb_{m_ik',m_in'}.
\]

\end{enumerate}
\end{thm}
\begin{proof}
By Lemma~\ref{lem:dregpt}, $N^\Eb_{k,n}$ is nonempty.
From Lemma~\ref{lem:regpt}, we have
\[
N^\Eb_{k,n}=\{R\in\Gc^\Eb_{k,n}\mid \frac nkR\mbox{ a regular point of }f\}.
\]
The regular value theorem now implies~(i).
Indeed,
\[
\Sc=\{P\in G^\Eb_{k,n}\mid P\mbox{ is a regular point of }f\}
\]
is an open subset and is therefore a regular, real analytic submanifold of the Grassman
manifold $G^\Eb_{k,n}$.
Now $c=(\frac nk)_{i=1}^k$ is a regular value of the restriction of $f$ to $\Sc$;
hence $N^\Eb_{k,n}$ is a regular, real analytic submanifold of $\Sc$ and thus also of $M_k(\Eb)$.

The assertions of~(ii) are clear with the possible exception of the inclusion $\subseteq$ in~\eqref{eq:Gcdecomp},
which we will prove by showing $R\in\Gc^\Eb_{k,n}$ implies $\sigma_R\in\Pc(k,k')$.
In fact, this is just a variant of the argument used to prove Theorem~\ref{thm:mfld}.
Let $P=\frac nkR$.
If $A\subseteq\{1,\ldots,k\}$ and $Q_AP=PQ_A$, then $PQ_A$ is a projection that
can be viewed as an $|A|\times|A|$ matrix, all of whose diagonal entries are $\frac nk=\frac{n'}{k'}$.
Hence the rank of $PQ_A$ is $|A|\frac{n'}{k'}$, which must therefore be an integer.
Since $\gcd(n',k')=1$, $|A|$ must be a multiple of $k'$.

For~(iii), we may without loss of generality assume the subsets $A_1,\ldots,A_\ell$ are
consequetive submintervals of $\{1,\ldots,k\}$.
Let $R\in\Gc^\Eb_{k,n}$.
Then $R\in N^\Eb_{k,n}$ if and only if $R$ is a block diagonal matrix $R=\diag(R_1,\ldots,R_\ell)$,
with $R_i\in N^\Eb_{m_ik',m_in'}$.
The assertions of~(iii) now follow readily from~(i).
\end{proof}

The above result on a manifold stratification structure of $\Gc^\Eb_{k,n}$,
together with the fiber bundle result Theorem~\ref{thm:fiber} (and Remark~\ref{rem:Cinf}),
yield directly a manifold stratification structure of $\Fc^\Eb_{k,n}$.
However, the following lemma allows us to describe this manifold stratification of $\Fc^\Eb_{k,n}$
directly in terms of frames.

\begin{lemma}\label{lem:frameinterp}
Let $F=(f_1,\ldots,f_k)\in\Fc^\Eb_{k,n}$ and let $A\subseteq\{1,\ldots,k\}$ be a subset.
Then $F^*F$ commutes with $Q_A$ if and only if there is a subspace $\Vc\subseteq\Eb^n$ such that
$(f_i)_{i\in A}$ forms a spherical tight frame for $\Vc$, while $(f_i)_{i\in A^c}$
forms a spherical tight frame for $\Vc^\perp$, where $A^c$ is the complement of $A$.
Moreover, if $F^*F$ commutes with $Q_A$, then the cardinality of $A$ is a multiple of $k/d$,
where $d=\gcd(k,n)$.
\end{lemma}
\begin{proof}
Multiplying $F$ on the right by a permutation matrix, if necessary,
we may without loss of generality assume $A=\{1,\ldots,p\}$ for some $p\in\{1,\ldots,k\}$.
Then $F=(F_1|F_2)$, where $F_1=(f_1,\ldots,f_p)$ and $F_2=(f_{p+1},\ldots,f_k)$,
and
\[
F^*F=\left(\begin{matrix}
F_1^*F_1 & F_1^*F_2 \\
F_2^*F_1 & F_2^*F_2
\end{matrix}\right).
\]
If $Q_A$ commutes with $F^*F$, then $F_1^*F_2=0$.
Moreover, since $\sqrt{\frac nk}F$ is a co--isometry, letting $\Vc$ be the range of $F_1^*$,
$F_1^*F_2=0$ implies that $\sqrt{\frac nk}F_1^*$ is an isometry from $\Eb^p$ onto $\Vc$,
while $\sqrt{\frac nk}F_2^*$ is an isometry from $\Eb^{k-p}$ onto $\Vc^\perp$,
i.e.\ $F_1$ is a spherical tight frame for $\Vc$ as is $F_2$ for $\Vc^\perp$.
The converse direction is clear.

An argument showing that $|A|$ must be a multiple of $k/d$ is contained in the proof of part~(ii)
of Theorem~\ref{thm:strat}.
\end{proof}

\begin{defi}\label{def:orthodecomp}\rm
Let $F=(f_i)_{i\in I}$ be a tight frame for some Hilbert space $\HEu$.
We say $F$ is {\em orthodecomposable} if there is a proper, nonempty subset $A\subseteq I$
such that $(f_i)_{i\in A}$ is a tight frame for some subspace $\Vc$ of $\HEu$, and $(f_i)_{i\in A^c}$
is a tight frame for $\Vc^\perp$.
\end{defi}

It is clear from Lemma~\ref{lem:frameinterp} and Definition~\ref{def:sigmaT}
(and, moreover, straighforward to show the analogous result directly in a more general context),
that for every $F\in\Fc^\Eb_{k,n}$, there is a unique
partition $\rho_F=\sigma_{F^*F}=\{A_1,\ldots,A_\ell\}$ of the set $\{1,\ldots,k\}$
and there is an orthogonal decomposition $\Eb^n=\Vc_1\oplus\cdots\oplus\Vc_\ell$
such that for every $j$, $(f_i)_{i\in A_j}$ is a STF for $\Vc_j$ that is not orthodecomposable.

We set
\[
\Mh^\Eb_{k,n}=\{F\in\Fc^\Eb_{k,n}\mid F\mbox{ not orthodecomposable }\},
\]
and, whenever $\sigma$ is a partition of the set $\{1,\ldots,k\}$, we set
\[
\Mh^\Eb_{k,n}(\sigma)=\{F\in\Fc^\Eb_{k,n}\mid\rho_F=\sigma\}.
\]

\begin{cor}\label{cor:Fstrat}
Let $k,n\in\Nats$ with $k>n$ and let $\Eb=\Reals$ or $\Eb=\Cpx$.
Let $m=n-1$ if $\Eb=\Reals$ and $m=2n-1$ if $\Eb=\Cpx$.
Then
\renewcommand{\labelenumi}{(\roman{enumi})}
\begin{enumerate}

\item
$\Mh^\Eb_{k,n}$ is a nonempty, regular, real analytic submanifold of $(S^m)^k$ with
\[
\dim(\Mh^\Eb_{k,n})=\begin{cases}
(k-\frac n2-1)(n-1),&\quad\Eb=\Reals, \\
2n(k-n)+n^2-k+1,&\quad\Eb=\Cpx.
\end{cases}
\]

\end{enumerate}
Let $d=\gcd(k,n)$ and let $k'=k/d$, $n'=n/d$.
Let $\Pc(k,k')$ be the set of all partitions of the set $\{1,\ldots,k\}$
into subsets whose cardinalities are multiples of $k'$.
Then
\renewcommand{\labelenumi}{(\roman{enumi})}
\begin{enumerate}
\stepcounter{enumi}

\item 
\[
\Fc^\Eb_{k,n}=\bigcup_{\sigma\in\Pc(k,k')}\Mh^\Eb_{k,n}(\sigma)
\]
and the sets $(\Mh^\Eb_{k,n}(\sigma))_{\sigma\in\Pc(k,k')}$ are pairwise disjoint;

\item
if $\sigma=\{A_1,\ldots,A_\ell\}\in\Pc(k,k')$ with $|A_i|=m_ik'$,
then $\Mh^\Eb_{k,n}(\sigma)$ is a nonempty, regular, real analytic submanifold of $(S^m)^k$
and is real--analytically diffeomorphic to the Cartesian product
\[
\prod_{i=1}^\ell\Mh^\Eb_{m_ik',m_in'}.
\]

\end{enumerate}
\end{cor}

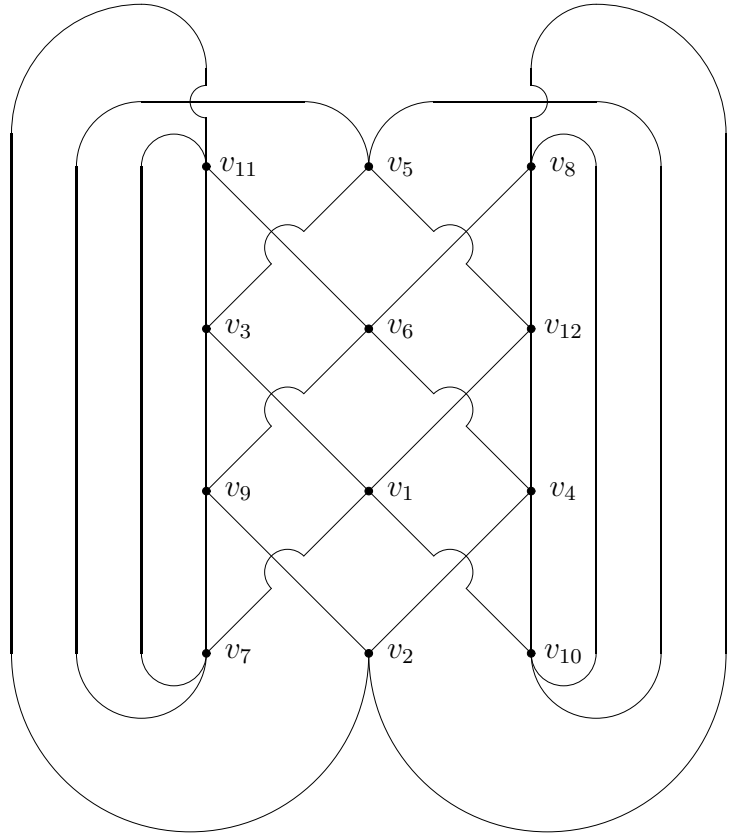
\section{The space $\Gc^\Reals_{4,2}$}\label{sec:G42}

In this section we will describe the space $\Gc^\Reals_{4,2}$ of
equivalence classes of spherical tight frames of four vectors in $\Reals^2$.
We begin with some general facts about spherical tight frames of $k$ vectors in $\Reals^2$.

The following proposition is elementary and well known; {\em cf}\/~\cite[Example 4.2]{BF} and~\cite[Thm 2.7]{GKK}. 
\begin{prop}\label{prop:R2C}
Let $k\in\Nats$, $k\ge2$ and let $f_1,\ldots,f_k\in\Reals^2$ with $f_1\ne0$.
Write $f_j=\left(\begin{smallmatrix}x_j \\ y_j\end{smallmatrix}\right)$ and let $z_j=x_j+iy_j\in\Cpx$.
Then $f_1,\ldots,f_k$ is a tight frame for $\Reals^2$ if and only if $\sum_{j=1}^kz_j^2=0$.
\end{prop}

\begin{cor}\label{cor:Fct}
Under the identification of $\Reals^2$ with $\Cpx$ used in Proposition~\ref{prop:R2C},
$\Fc^\Reals_{k,2}$ is identified with
\[
\Fct_{k,2}=\{(z_1,\ldots,z_k)\in\Tcirc^k\mid\sum_{j=1}^kz_j^2=0\}.
\]
Furthermore, the orbit space $\Gc^\Reals_{k,2}=\Fc^\Reals_{k,2}/\Oc^\Reals_2$ is identified with
the orbit space $\Gct_{k,2}$ of $\Fct_{k,2}$ under the group of transformations
generated the action of $\Ints_2$ by complex conjugation,
\[
(z_1,\ldots,z_k)\mapsto(\overline{z_1},\ldots,\overline{z_k})
\]
and the action of $\Tcirc$ by rotations,
\[
e^{i\theta}\cdot(z_1,\ldots,z_k)=(e^{i\theta}z_1,\ldots,e^{i\theta}z_k).
\]
\end{cor}

\begin{figure}[bp]
\caption{The image of $\tau_{p,\eps}$ as a graph.}\label{fig:taugraph}
\setlength{\unitlength}{0.92in}
\begin{picture}(3.5,3)(-0.75,-0.5)
\put(1,0){\circle*{0.05}}
\put(0,1){\circle*{0.05}}
\put(2,1){\circle*{0.05}}
\put(1,2){\circle*{0.05}}
\put(1,0){\line(1,1){1}}
\put(1,0){\line(-1,1){1}}
\put(1,2){\line(1,-1){1}}
\put(1,2){\line(-1,-1){1}}
\put(-0.7,1){$\tau_{p,\eps}(-1)$}
\put(2.1,1){$\tau_{p,\eps}(1)$}
\put(1.1,2.1){$\tau_{p,\eps}(i)$}
\put(1.1,-0.15){$\tau_{p,\eps}(-i)$}
\end{picture}
\end{figure}
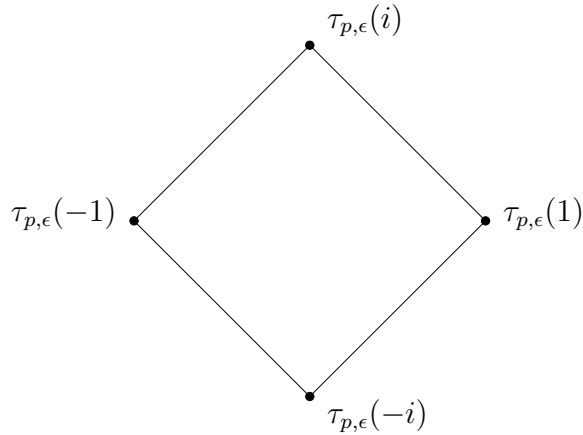

Since any four elements of $\Tcirc$ that sum to zero can be divided into pairs that
are negatives of each other, from Corollary~\ref{cor:Fct} we easily prove the well--known
result that any spherical tight frame of four vectors in $\Reals^2$ consists of two
orthonormal bases.
Thus, the orbit space $\Fct_{4,2}/\Tcirc$ is the union of the images of the twelve maps
\[
\big(\sigma_{p,\eps_1,\eps_2}\big)_{2\le p\le 4,\,\eps_1,\eps_2\in\{\pm1\},}
\]
where $\sigma_{p,\eps_1,\eps_2}:\Tcirc\to\Fct_{4,2}/\Tcirc$ are given by
\begin{align*}
\sigma_{2,\eps_1,\eps_2}(\zeta)&=[(1,\eps_1i,\zeta,\eps_2i\zeta)] \\
\sigma_{3,\eps_1,\eps_2}(\zeta)&=[(1,\zeta,\eps_1i,\eps_2i\zeta)] \\
\sigma_{4,\eps_1,\eps_2}(\zeta)&=[(1,\zeta,\eps_2i\zeta,\eps_1i)].
\end{align*}
For every $p$ we have
\[
\overline{\sigma_{p,\eps_1,\eps_2}(\zeta)}=\sigma_{p,-\eps_1,-\eps_2}(\overline\zeta).
\]
Thus $\Gct_{4,2}=(\Fct_{4,2}/\Tcirc)/\Ints_2$ is the union of the images of the
six maps
\begin{equation}\label{eq:taus}
\big(\tau_{p,\eps}\big)_{2\le p\le 4,\,\eps\in\{\pm1\},}
\end{equation}
where $\tau_{p,\eps}:\Tcirc\to\Gct_{4,2}$ and $\tau_{p,\eps}(\zeta)=[\mu_{p,\eps}(\zeta)]$
with
\begin{align*}
\mu_{2,\eps}(\zeta)&=(1,i,\zeta,\eps i\zeta) \\
\mu_{3,\eps}(\zeta)&=(1,\zeta,i,\eps i\zeta) \\
\mu_{4,\eps}(\zeta)&=(1,\zeta,\eps i\zeta,i).
\end{align*}
Note we have $\tau_{p,\eps}(\zeta)=\tau_{p',\eps'}(\zeta')$ if and only if either
$\mu_{p,\eps}(\zeta)=\mu_{p',\eps'}(\zeta')$ or $\mu_{p,\eps}(\zeta)=\overline{\mu_{p',\eps'}(\zeta')}$.
For brevity, we will abreviate $\eps=\pm1$ by $\eps=\pm$.
For each $p$ and $\eps$, we have
\begin{alignat*}{2}
\mu_{p,\eps}(\zeta)&=\mu_{p,\eps}(\zeta')\quad&\Leftrightarrow\quad\zeta=\zeta' \\[1ex]
\mu_{p,\eps}(\zeta)&=\overline{\mu_{p,\eps}(\zeta')}\quad&\mbox{never happens.}
\end{alignat*}
Thus each $\tau_{p,\eps}$ is injective.
It follows that $\Gct_{4,2}$ is the indentification space of six circles, glued together
according to how the images of the six maps~\eqref{eq:taus} overlap.
For each $p$,
\begin{alignat*}{2}
\mu_{p,+}(\zeta)&=\mu_{p,-}(\zeta')\quad&&\mbox{never happens} \\[1ex]
\mu_{p,+}(\zeta)&=\overline{\mu_{p,-}(\zeta')}\quad&&\mbox{never happens}
\end{alignat*}
so $\tau_{p,+}$ and $\tau_{p,-}$ have disjoint images.
Furthermore,
\begin{alignat*}{2}
\mu_{2,\eps}(\zeta)&=\mu_{3,\eps'}(\zeta')\quad
 &&\Leftrightarrow\quad\zeta=\zeta'=i,\,\eps=\eps' \\[1ex]
\mu_{2,\eps}(\zeta)&=\overline{\mu_{3,\eps'}(\zeta')}\quad
 &&\Leftrightarrow\quad\zeta=\zeta'=-i,\,\eps=\eps' \\[1ex]
\mu_{2,\eps}(\zeta)&=\mu_{4,\eps'}(\zeta')\quad
 &&\Leftrightarrow\quad\zeta=-\eps',\,\zeta'=i,\,\eps=-\eps' \\[1ex]
\mu_{2,\eps}(\zeta)&=\overline{\mu_{4,\eps'}(\zeta')}\quad
 &&\Leftrightarrow\quad\zeta=\eps',\,\zeta'=-i,\,\eps=-\eps' \\[1ex]
\mu_{3,\eps}(\zeta)&=\mu_{4,\eps'}(\zeta')\quad
 &&\Leftrightarrow\quad\zeta=\zeta'=\eps=\eps' \\[1ex]
\mu_{3,\eps}(\zeta)&=\overline{\mu_{4,\eps'}(\zeta')}\quad
 &&\Leftrightarrow\quad\zeta=\zeta'=-\eps,\,\eps=\eps'.
\end{alignat*}
Representing the image of $\tau_{p,\eps}$ as a graph in the obvious way depicted in
Figure~\ref{fig:taugraph}, the only identifications that occur in $\Gct_{4,2}$
are among the labeled vertices.
In particular,
we find that $\Gct_{4,2}$ is homeomorphic to the quotient graph whose vertices are
\begin{alignat*}{3}
v_1&=\{\tau_{2,+}(1),\tau_{4,-}(i)\}\quad
&v_2&=\{\tau_{2,-}(1),\tau_{4,+}(-i)\}\quad
&v_3&=\{\tau_{2,+}(i),\tau_{3,+}(i)\} \\
v_4&=\{\tau_{2,-}(i),\tau_{3,-}(i)\}\quad
&v_5&=\{\tau_{2,+}(-1),\tau_{4,-}(-i)\}\quad
&v_6&=\{\tau_{2,-}(-1),\tau_{4,+}(i)\} \\
v_7&=\{\tau_{2,+}(-i),\tau_{3,+}(-i)\}\quad
&v_8&=\{\tau_{2,-}(-i),\tau_{3,-}(-i)\}\quad
&v_9&=\{\tau_{3,+}(1),\tau_{4,+}(1)\} \\
v_{10}&=\{\tau_{3,-}(1),\tau_{4,-}(1)\}\quad
&v_{11}&=\{\tau_{3,+}(-1),\tau_{4,+}(-1)\}\quad
&v_{12}&=\{\tau_{3,-}(-1),\tau_{4,-}(-1)\},
\end{alignat*}
and where each vertex has four edges, surviving from the original graphs in Figure~\ref{fig:taugraph}.
In conclusion:
\begin{thm}\label{thm:G42}
the space $\Gc^\Reals_{4,2}$ is homeomorphic to the graph with twelve vertices and twenty--four
edges that is depicted in Figure~\ref{fig:G42graph}.
\begin{figure}[bp]
\caption{The space $\Gc^\Reals_{4,2}$.}\label{fig:G42graph}
\setlength{\unitlength}{0.85in}
\begin{picture}(4,5.5)(-1,-1.3)
\multiput(0,0)(0,1){4}{\circle*{0.05}}
\multiputlist(0.2,0)(0,1){$v_7$,$v_9$,$v_3$,$v_{11}$}
\multiput(0,0)(0,1){3}{\line(0,1){1}}
\multiput(0,0)(0,1){3}{\line(1,1){0.4}}
\multiput(0.6,0.6)(0,1){3}{\line(1,1){0.4}}
\multiput(0.5,0.5)(0,1){3}{\arc{0.28284}{2.35619}{5.49779}}
\multiput(1,0)(0,1){4}{\circle*{0.05}}
\multiputlist(1.2,0)(0,1){$v_2$,$v_1$,$v_6$,$v_5$}
\multiput(1,0)(0,1){3}{\line(-1,1){1}}
\multiput(1,0)(0,1){3}{\line(1,1){1}}
\multiput(2,0)(0,1){4}{\circle*{0.05}}
\multiputlist(2.2,0)(0,1){$v_{10}$,$v_4$,$v_{12}$,$v_8$}
\multiput(2,0)(0,1){3}{\line(0,1){1}}
\multiput(2,0)(0,1){3}{\line(-1,1){0.4}}
\multiput(1.4,0.6)(0,1){3}{\line(-1,1){0.4}}
\multiput(1.5,0.5)(0,1){3}{\arc{0.28284}{3.92699}{7.06858}}
\put(-0.2,3){\arc{0.4}{3.14159}{6.28319}}
\put(-0.4,0){\line(0,1){3}}
\put(-0.2,0){\arc{0.4}{0}{3.14159}}
\put(0.6,3){\arc{0.8}{4.71239}{0}}
\put(0.6,3.4){\line(-1,0){1}}
\put(-0.4,3){\arc{0.8}{3.14159}{4.71239}}
\put(-0.8,0){\line(0,1){3}}
\put(-0.4,0){\arc{0.8}{0}{3.14159}}
\put(0,3){\line(0,1){0.3}}
\put(0,3.4){\arc{0.2}{1.57080}{4.71239}}
\put(0,3.5){\line(0,1){0.1}}
\put(-0.4,3.6){\arc{0.8}{4.71239}{0}}
\put(-0.4,3.2){\arc{1.6}{3.14159}{4.71239}}
\put(-1.2,0){\line(0,1){3.2}}
\put(-0.1,0){\arc{2.2}{0}{3.14159}}
\put(2.2,3){\arc{0.4}{3.14159}{6.28319}}
\put(2.4,0){\line(0,1){3}}
\put(2.2,0){\arc{0.4}{0}{3.14159}}
\put(1.4,3){\arc{0.8}{3.14159}{4.71239}}
\put(1.4,3.4){\line(1,0){1}}
\put(2.4,3){\arc{0.8}{4.71239}{6.28319}}
\put(2.8,0){\line(0,1){3}}
\put(2.4,0){\arc{0.8}{0}{3.14159}}
\put(2,3){\line(0,1){0.3}}
\put(2,3.4){\arc{0.2}{4.71239}{7.85398}}
\put(2,3.5){\line(0,1){0.1}}
\put(2.4,3.6){\arc{0.8}{3.14159}{4.71239}}
\put(2.4,3.2){\arc{1.6}{4.71239}{6.28319}}
\put(3.2,0){\line(0,1){3.2}}
\put(2.1,0){\arc{2.2}{0}{3.14159}}
\end{picture}
\end{figure}
\end{thm}

\section{The space $\Gc^\Reals_{5,2}$}\label{sec:G52}

In this section we will describe the space $\Gc^\Reals_{5,2}$ of
equivalence classes of tight spherical frames of five vectors in $\Reals^2$.
By Corollary~\ref{cor:Fct}, $\Fc^\Reals_{5,2}$ is homeomorphic to
$\Fct_{5,2}$.
We have
\[
\Ec\eqdef\Fct_{5,2}/\Tcirc=\{(z_1,z_2,z_3,z_4)\in\Tcirc^4\mid\sum_{j=1}^4z_j^2=-1\}
\]
and $\Gc^\Reals_{5,2}$ is homeomorphic to $\Ec/\Ints_2$ where $\Ints_2$ acts on $\Ec$ by complex conjugation.
Let
\[
\Dc=\{(w_1,w_2,w_3,w_4)\in\Tcirc^4\mid\sum_{j=1}^4w_j=-1\}
\]
and let $p:\Ec\to\Dc$ be $p(z_1,z_2,z_3,z_4)=(z_1^2,z_2^2,z_3^2,z_4^2)$.
Then $p$ is a sixteen--fold covering map and $p$ intertwines complex conjugation with complex conjugation.

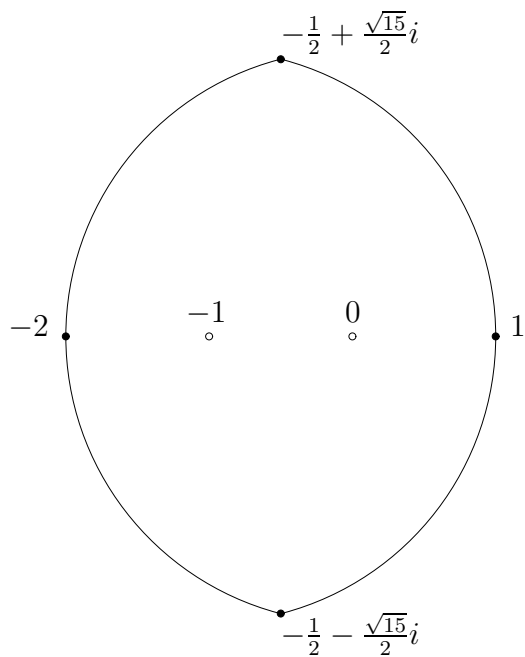
\begin{figure}[bp]
\caption{The set $A$.}\label{fig:A}
\setlength{\unitlength}{0.75in}
\begin{picture}(5,5)(-3,-2.5)
\put(0,0){\arc{4}{1.83220}{4.45098}}
\put(-1,0){\arc{4}{4.97380}{7.59258}}
\put(-1,0){\circle{0.05}}
\put(-1.16,0.1){$-1$}
\put(0,0){\circle{0.05}}
\put(-0.05,0.1){$0$}
\put(-2,0){\circle*{0.05}}
\put(-2.4,0){$-2$}
\put(1,0){\circle*{0.05}}
\put(1.1,0){$1$}
\put(-0.5,1.936){\circle*{0.05}}
\put(-0.5,2.05){$-\frac12+\frac{\sqrt{15}}2i$}
\put(-0.5,-1.936){\circle*{0.05}}
\put(-0.5,-2.15){$-\frac12-\frac{\sqrt{15}}2i$}
\end{picture}
\end{figure}

\begin{figure}[bp]
\caption{The space $B$.}\label{fig:B}
\setlength{\unitlength}{0.75in}
\begin{picture}(5,5)(-3,-2.5)
\put(0,0){\arc{4}{1.83220}{4.45098}}
\put(-2.4,0){$\partial_\ell B$}
\put(-1,0){\arc{4}{4.97380}{7.59258}}
\put(1.1,0){$\partial_r B$}
\put(-1,0){\circle{0.7}}
\put(-1.1,0.4){$C_{-1}$}
\put(0,0){\circle{0.7}}
\put(-0.05,0.4){$C_0$}
\put(-0.5,1.936){\circle*{0.05}}
\put(-0.5,2.05){$q$}
\put(-0.5,-1.936){\circle*{0.05}}
\put(-0.5,-2.15){$\overline q$}
\end{picture}
\end{figure}
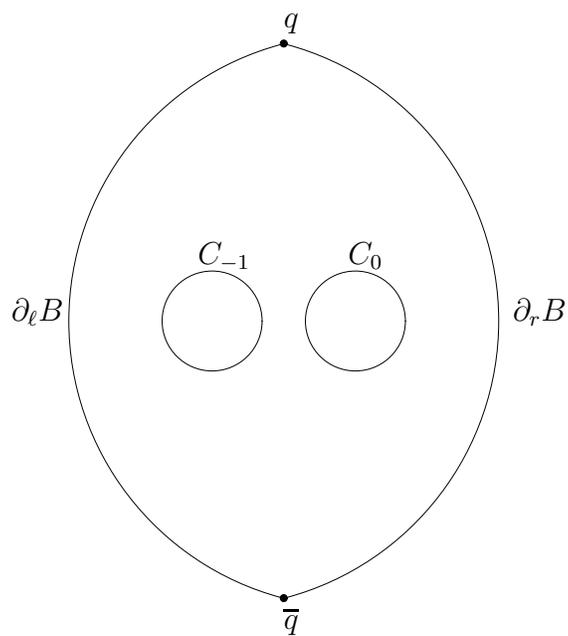

Let
\[
A=\{a\in\Cpx\mid0<|a|\le2,\,0<|-1-a|\le2\}.
\]
The set $A$ is pictured in Figure~\ref{fig:A},
where of course we have $-1,0\notin A$.
Let $B$ be the topological space obtained from $A$ by dilating the punctures at $-1$ and $0$
and gluing copies $C_{-1}$, respectively $C_0$, of the circle onto the boundaries of the resulting holes.
The space $B$ is pictured in Figure~\ref{fig:B}.
More formally, as a set we let $B$ be the disjoint union of the ranges of three
injective maps,
\[
\alpha:A\to B,\qquad
\tau_{-1}:\Tcirc\to B,\qquad
\tau_0:\Tcirc\to B,
\]
where $C_j$ is the image of $\tau_j$,
and where the topology of $B$ is defined as follows:
\renewcommand{\labelenumi}{$\bullet$}
\begin{enumerate}

\item
a neighborhood of $\alpha(a)$ in $B$ is a subset of $B$ containing
\[
\alpha(\{a'\in A\mid\,|a'-a|<\eps\})
\]
for some $\eps>0$;

\item
a neighborhood of $\tau_0(e^{i\theta})$ in $B$ is a subset of $B$ containing
\begin{align*}
\{\tau_0(e^{i\theta'})&\mid\theta'\in\Reals,\,|\theta'-\theta|<\eps\}\cup \\
&\cup\{\alpha(re^{i\theta'})\mid\theta'\in\Reals,\,|\theta'-\theta|<\eps,\,0<r<\eps\}
\end{align*}
for some $0<\eps<1$;

\item
a neighborhood of $\tau_{-1}(e^{i\theta})$ in $B$ is a subset of $B$ containing
\begin{align*}
\{\tau_{-1}(e^{i\theta'})&\mid\theta'\in\Reals,\,|\theta'-\theta|<\eps\}\cup \\
&\cup\{\alpha(-1+re^{i\theta'})\mid\theta'\in\Reals,\,|\theta'-\theta|<\eps,\,0<r<\eps\}
\end{align*}
for some $0<\eps<1$.
\end{enumerate}
In Figure~\ref{fig:B}, we have
$q=\alpha(\frac12+\frac{\sqrt{15}}2i)$ and $\overline q=\alpha(\frac12-\frac{\sqrt{15}}2i)$.
Then $B$ is a compact Hausdorff space.
Let $p_A:B\to\cl(A)$, where $\cl(A)=A\cup\{0,1\}$ denotes the closure of $A$,
be the continuous map given by
\[
p_A(\alpha(a))=a,\qquad p_A(\tau_{-1}(\zeta))=-1,\qquad p_A(\tau_0(\zeta))=0,\qquad(\zeta\in\Tcirc).
\]

Given $a\in A$, there are precisely two values of $w\in\Tcirc$ such that there exists $w'\in\Tcirc$
with $w+w'=a$.
These are illustrated in Figure~\ref{fig:w}.
\begin{figure}[bp]
\caption{$w+w'=a$.}\label{fig:w}
\setlength{\unitlength}{1in}
\begin{picture}(2,2.5)(-1.5,-1)
\put(0,0){\circle*{0.05}}
\put(0.07,-0.07){$0$}
\put(-0.577,0.333){\circle*{0.05}}
\put(0,0){\line(-577,333){0.577}}
\put(-0.70,0.35){$a$}
\put(0.183,0.983){\circle*{0.05}}
\put(0,0){\line(183,983){0.183}}
\put(0.183,0.983){\line(-76,-65){0.76}}
\put(0.1,1.05){$w$}
\put(-0.76,-0.65){\circle*{0.05}}
\put(0,0){\line(-76,-65){0.76}}
\put(-0.76,-0.65){\line(183,983){0.183}}
\put(-1.0,-0.75){$w'$}
\end{picture}
\end{figure}
Let $w_\ell(a)$ be the one of these two values that lies to the left as one travels the line segment
from $0$ to $a$, and let $w_r(a)$ be the other value, lying to the right.
In Figure~\ref{fig:w}, $w_\ell(a)$ is the point labeled $w'$ and $w_r(a)$ is the point labeled $w$.
Note we have
\begin{equation}\label{eq:wconj}
w_\ell(\overline a)=\overline{w_r(a)},\qquad w_r(\overline a)=\overline{w_\ell(a)}.
\end{equation}

Let $\phi:B\to\Dc$ be given by
\begin{align}
\phi(\alpha(a))&=(w_\ell(a),w_r(a),w_\ell(-1-a),w_r(-1-a)),\qquad(a\in A) \label{eq:phia} \\
\phi(\tau_0(\zeta))&=(i\zeta,-i\zeta,-\tfrac12-\tfrac{\sqrt3}2i,-\tfrac12+\tfrac{\sqrt3}2i)
 \qquad(\zeta\in\Tcirc). \label{eq:phi0} \\
\phi(\tau_{-1}(\zeta))&=(-\tfrac12-\tfrac{\sqrt3}2i,-\tfrac12+\tfrac{\sqrt3}2i,-i\zeta,i\zeta)
 \qquad(\zeta\in\Tcirc) \label{eq:phi-1}
\end{align}
Then $\phi$ is injective and continuous, hence a homeomorphism onto its image.
Let $V=\{e,t,u,v\}$ be the Klein 4--group with $e$ the identity element and let $V$ act on $\Tcirc^4$
by
\begin{align*}
t\odot(\zeta_1,\zeta_2,\zeta_3,\zeta_4)=(\zeta_2,\zeta_1,\zeta_3,\zeta_4) \\
u\odot(\zeta_1,\zeta_2,\zeta_3,\zeta_4)=(\zeta_1,\zeta_2,\zeta_4,\zeta_3) \\
v\odot(\zeta_1,\zeta_2,\zeta_3,\zeta_4)=(\zeta_2,\zeta_1,\zeta_4,\zeta_3).
\end{align*}
For $g\in V$, let $\phi_g:B\to\Dc$ be $\phi_g(b)=g\odot(\phi(b))$.
Then
\begin{equation}\label{eq:D}
\Dc=\bigcup_{g\in V}\phi_g(B).
\end{equation}
Let
\begin{align*}
\partial_\ell B&=\{\alpha(a)\mid a\in A,\,|a|=2\} \\
\partial_r B&=\{\alpha(a)\mid a\in A,\,|-1-a|=2\}
\end{align*}
be the left and right boundaries of $B$, as indicated in Figure~\ref{fig:B}.
Thus $\partial_\ell B\cap\partial_r B=\{q,\overline q\}$.

\begin{lemma}\label{lem:phiB}
Let $b,b'\in B$, $g,g'\in V$.
Then
\begin{equation}\label{eq:phigb}
\phi_g(b)=\phi_{g'}(b')
\end{equation}
if and only if at least one of the following holds:
\renewcommand{\labelenumi}{(\roman{enumi})}
\begin{enumerate}

\item $b=b'$ and $g=g'$

\item $b=b'\in\partial_\ell B\cap\partial_r B$

\item $b=b'\in\partial_\ell B$ and $g^{-1}g'=t$

\item $b=b'\in\partial_r B$ and $g^{-1}g'=u$

\item $b=\tau_0(\zeta)$ and $b'=\tau_0(-\zeta)$ for some $\zeta\in\Tcirc$, and $g^{-1}g'=t$

\item $b=\tau_{-1}(\zeta)$ and $b'=\tau_{-1}(-\zeta)$ for some $\zeta\in\Tcirc$, and $g^{-1}g'=u$.

\end{enumerate}
\end{lemma}
\begin{proof}
We write $\phi(b)=(w_1(b),w_2(b),w_3(b),w_4(b))$ for $b\in B$.
Then $w_1(b)+w_2(b)=p_A(b)$, so~\eqref{eq:phigb} implies $p_A(b)=p_A(b')$,
and one of the following holds:
\begin{gather}
b=b'=\alpha(a),\quad\mbox{some }a\in A \label{eq:ba} \\
b,b'\in C_0 \label{eq:bC0} \\
b,b'\in C_{-1}. \label{eq:bC-1}
\end{gather}
Suppose~\eqref{eq:ba} holds.
Since $w_1(b)=w_2(b)$ if and only if $b\in\partial_\ell(B)$
and $w_3(b)=w_4(b)$ if and only if $b\in\partial_r(B)$, we quickly deduce that
at least one of conditions (i)--(iv) holds,
and conversely, any of (i)--(iv) implies~\eqref{eq:phigb}.
If~\eqref{eq:bC0} holds, then using~\eqref{eq:phi0} we get \eqref{eq:phigb}$\Leftrightarrow$(i) or (v),
and if~\eqref{eq:bC-1} holds, then using~\eqref{eq:phi-1} we get \eqref{eq:phigb}$\Leftrightarrow$(i) or (vi).
\end{proof}

By~\eqref{eq:D} and the fact that each $\phi_g$ is injective,
the space $\Dc$ is obtained by gluing together four copies of $B$ according to how their
images under $\phi_e,\phi_t,\phi_u$ and $\phi_v$ overlap, as described in Lemma~\ref{lem:phiB}.
We will not pursue this, but we will use a similar reasoning to investigate the space $\Ec/\Ints_2$.

Let $\beta:B\to B$ be the homeomorphism of order two defined by complex conjugation, namely
\begin{alignat*}{2}
\beta(\alpha(a))&=\alpha(\overline a)&\quad&(a\in A) \\
\beta(\tau_{-1}(\zeta))&=\tau_{-1}(\overline\zeta)&&(\zeta\in\Tcirc) \\
\beta(\tau_0(\zeta))&=\tau_0(\overline\zeta)&&(\zeta\in\Tcirc).
\end{alignat*}
Then $\beta^2=\id$.
If $b\in B$ and $\phi(b)=(w_1,w_2,w_3,w_4)$, then from~\eqref{eq:wconj}--\eqref{eq:phi0} we see
\begin{equation}\label{eq:phibeta}
\phi(\beta(b))=(\overline{w_2},\overline{w_1},\overline{w_4},\overline{w_3})=\overline{\phi_v(b)}.
\end{equation}
Let $\Bt$ be the upper half of $B$, namely
\begin{align*}
\Bt={}&\{\alpha(a)\mid a\in A,\,\RealPart a\ge0\}\cup \\
&\cup\{\tau_{-1}(\zeta)\mid\zeta\in\Tcirc,\,\RealPart\zeta\ge0\}\cup \\
&\cup\{\tau_0(\zeta)\mid\zeta\in\Tcirc,\,\RealPart\zeta\ge0\}.
\end{align*}
Then $B=\Bt\cup\beta(\Bt)$.

For future use, we want to consider also the homeomorphism $\gamma:B\to B$ given by
\begin{alignat*}{2}
\gamma(\alpha(a))&=\alpha(-1-\overline a)&\quad&(a\in A) \\
\gamma(\tau_0(\zeta))&=\tau_{-1}(-\overline\zeta)&&(\zeta\in\Tcirc) \\
\gamma(\tau_{-1}(\zeta))&=\tau_0(-\overline\zeta)&&(\zeta\in\Tcirc).
\end{alignat*}
Then $\gamma^2=\id$.
Moreover, for $b\in B$, if $\phi(b)=(w_1,w_2,w_3,w_4)$, then
\begin{equation}\label{eq:phigamma}
\phi(\gamma(b))=(\overline{w_4},\overline{w_3},\overline{w_2},\overline{w_1}).
\end{equation}

We have
\[
\phi(q)=(-\tfrac14+\tfrac{\sqrt{15}}4i,-\tfrac14+\tfrac{\sqrt{15}}4i,
-\tfrac14-\tfrac{\sqrt{15}}4i,-\tfrac14-\tfrac{\sqrt{15}}4i).
\]
By homotopy lifting, there is a unique continuous map $\phit:\Bt\to\Ec$
such that $p\circ\phit=\phi\restrict_\Bt$ and such that
\[
\phit(q)=(\sqrt{\tfrac38}+\sqrt{\tfrac58}\,i,\sqrt{\tfrac38}+\sqrt{\tfrac58}\,i,
\sqrt{\tfrac38}-\sqrt{\tfrac58}\,i,\sqrt{\tfrac38}-\sqrt{\tfrac58}\,i).
\]
Let $E$ be the multiplicative subgroup
\[
E=\{(\eps_1,\eps_2,\eps_3,\eps_4)\mid\eps_j\in\{\pm1\}\}
\]
of $\Tcirc^4$.
Given $g\in V$ and $\eps\in E$, let $\phit_g^\eps:\Bt\to\Ec$ be
\[
\phit_g^\eps(b)=\eps\cdot(g\odot(\phit(b))).
\]
Then for any $g\in V$, $(\phit_g^\eps)_{\eps\in E}$ are the sixteen different liftings
of $\phi_g\restrict_\Bt$ to $\Ec$ (under the covering projection $p$).
in particular, for $b\in\Bt$,
\begin{equation}\label{eq:philift}
p^{-1}(\phi_g(b))=\{\phit_g^\eps(b)\mid\eps\in E\}.
\end{equation}
If $b\in\beta(\Bt)$, then by~\eqref{eq:phibeta},
\[
p(\overline{\phit_g^\eps(\beta(b))})=\overline{p(\phit_g^\eps(\beta(b)))}
=\overline{\phi_g(\beta(b))}=\phi_{vg}(b).
\]
Therefore,
\begin{equation}\label{eq:phibetalift}
p^{-1}(\phi_{vg}(b))=\{\overline{\phit_g^\eps(\beta(b))}\mid\eps\in E\}.
\end{equation}
Let $\psi_g^\eps:\Bt\to\Ec/\Ints_2$ be $\phit_g^\eps$ followed by the quotient map $\Ec\to\Ec/\Ints_2$
of $\Ec$ under the action of complex conjugation.
Because $B=\Bt\cup\beta(\Bt)$, from~\eqref{eq:philift} and~\eqref{eq:phibetalift}, we get
that $\Ec/\Ints_2$ is covered by the images of the sixty--four maps
$(\psi_g^\eps)_{g\in V,\,\eps\in E}$.
Our goal is to understand how these images overlap.

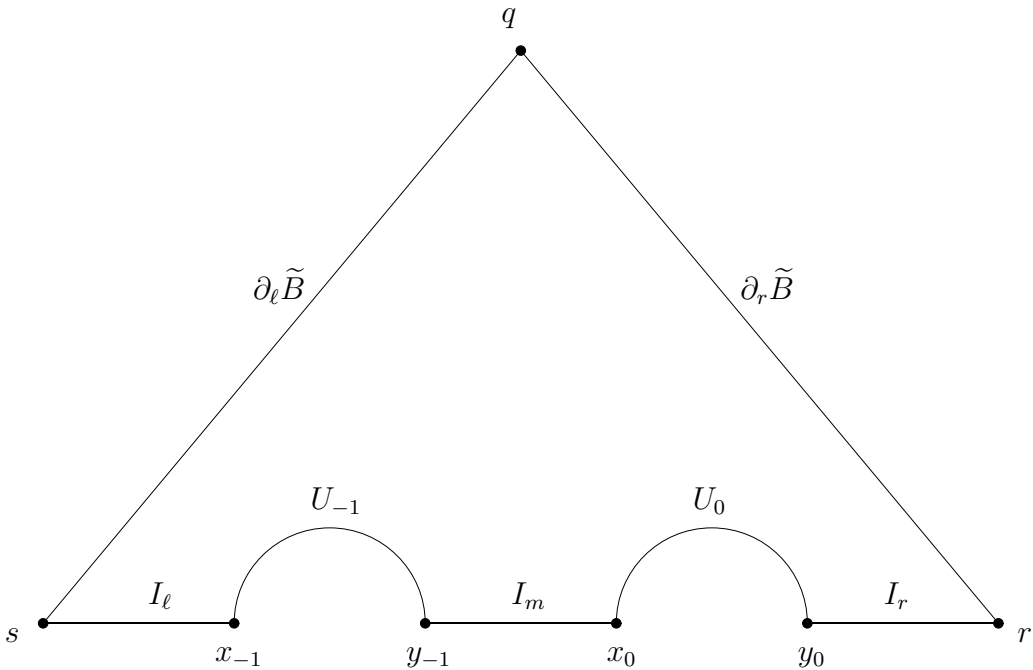
\begin{figure}[bp]
\caption{The space $\Bt$ with labeled boundary.}\label{fig:Bt}
\setlength{\unitlength}{1in}
\begin{picture}(5.5,3.5)(-2.75,0)
\put(0,3){\circle*{0.05}}
\put(-0.1,3.15){$q$}
\put(0,3){\line(-5,-6){2.5}}
\put(-1.4,1.7){$\partial_\ell\Bt$}
\put(0,3){\line(5,-6){2.5}}
\put(1.15,1.7){$\partial_r\Bt$}
\put(-2.5,0){\circle*{0.05}}
\put(-2.7,-0.1){$s$}
\put(-2.5,0){\line(1,0){1}}
\put(-1.95,0.1){$I_\ell$}
\put(-1.5,0){\circle*{0.05}}
\put(-1.6,-0.2){$x_{-1}$}
\put(-1,0){\arc{1}{3.14159}{6.28319}}
\put(-1.1,0.6){$U_{-1}$}
\put(-0.5,0){\circle*{0.05}}
\put(-0.6,-0.2){$y_{-1}$}
\put(-0.5,0){\line(1,0){1}}
\put(-0.05,0.1){$I_m$}
\put(0.5,0){\circle*{0.05}}
\put(0.45,-0.2){$x_0$}
\put(1,0){\arc{1}{3.14159}{6.28319}}
\put(0.9,0.6){$U_0$}
\put(1.5,0){\circle*{0.05}}
\put(1.45,-0.2){$y_0$}
\put(1.5,0){\line(1,0){1}}
\put(1.9,0.1){$I_r$}
\put(2.5,0){\circle*{0.05}}
\put(2.6,-0.1){$r$}
\end{picture}
\end{figure}

We shall compute the values of $\phit$ on the boundary of $\Bt$.
The space $\Bt$ is pictured and its boundary labeled in Figure~\ref{fig:Bt}.
In particular we have the indicated points
\begin{gather*}
q=\alpha(-\tfrac12+\tfrac{\sqrt{15}}2i) \\
s=\alpha(-2),\quad r=\alpha(1) \\
x_j=\tau_j(-1),\quad y_j=\tau_j(1),\quad(j\in\{-1,0\})
\end{gather*}
and the (closed) intervals and arcs
\begin{align*}
\partial_z\Bt&=\partial_zB\cap\Bt,\quad(z=\ell,r) \\
U_j&=C_j\cap\Bt,\quad(j\in\{-1,0\}) \\
I_\ell&=\alpha([-2,-1))\cup\{x_{-1}\} \\
I_m&=\{y_{-1}\}\cup\alpha((-1,0))\cup\{x_0\} \\
I_r&=\{y_0\}\cup\alpha((0,1]).
\end{align*}
Let us first determine $\phit(\alpha(-\frac12+it))$ for $t\in[0,\frac{\sqrt{15}}2]$.
Let $b=\alpha(-\frac12+it)$ and write
\[
\phi(b)=(w_1(b),w_2(b),w_3(b),w_4(b)),\qquad\phit(b)=(z_1(b),z_2(b),z_3(b),z_4(b)).
\]
As $t$ descends from $\frac{\sqrt{15}}2$ to $0$, $w_1(b)$ moves from $-\frac14+\frac{\sqrt{15}}4i$
to $-\frac14-\frac{\sqrt{15}}4i$, avoiding the first quadrant, and $w_2(b)$ moves from
$-\frac14+\frac{\sqrt{15}}4i$ back to $-\frac14+\frac{\sqrt{15}}4i$, staying in the upper half--plane.
Therefore, $z_1(b)$ changes from $\sqrt{\frac38}+\sqrt{\frac58}\,i$ to $-\sqrt{\frac38}+\sqrt{\frac58}\,i$
and $z_2(b)$ moves from $\sqrt{\frac38}+\sqrt{\frac58}\,i$ back to $\sqrt{\frac38}+\sqrt{\frac58}\,i$.
Finally, since $b=\gamma(b)$, from~\eqref{eq:phigamma} we have $w_3(b)=\overline{w_2(b)}$ and
$w_4(b)=\overline{w_1(b)}$; consequently  $z_3(b)=\overline{z_2(b)}$ and
$z_4(b)=\overline{z_1(b)}$.
Hence
\[
\phit(\alpha(-\tfrac12))=\Big(-\sqrt{\tfrac38}+\sqrt{\tfrac58}\,i,\,\sqrt{\tfrac38}+\sqrt{\tfrac58}\,i,\,
\sqrt{\tfrac38}-\sqrt{\tfrac58}\,i,\,-\sqrt{\tfrac38}-\sqrt{\tfrac58}\,i\Big).
\]
Now it is easy to calculate the values of $\phi$ and $\phit$ on the bottom part of the boundary
of $\Bt$.
On $I_m$, for $t\in(-1,0)$,
\begin{align}
\phi(\alpha(t))&=\Big(\tfrac t2-\tfrac{\sqrt{4-t^2}}2\,i,\,\tfrac t2+\tfrac{\sqrt{4-t^2}}2\,i,\,
\tfrac{-1-t}2-\tfrac{\sqrt{4-(1+t)^2}}2\,i,\,\tfrac{-1-t}2+\tfrac{\sqrt{4-(1+t)^2}}2\,i\Big) \notag \\[0.5ex]
\phit(\alpha(t))&=\begin{aligned}[t]\Big(&-\sqrt{\tfrac12+\tfrac t4}+\sqrt{\tfrac12-\tfrac t4}\,i,\,
\sqrt{\tfrac12+\tfrac t4}+\sqrt{\tfrac12-\tfrac t4}\,i,\, \\
&\sqrt{\tfrac12-\tfrac{t+1}4}-\sqrt{\tfrac12+\tfrac{t+1}4}\,i,\,
-\sqrt{\tfrac12-\tfrac{t+1}4}-\sqrt{\tfrac12+\tfrac{t+1}4}\,i\Big).\end{aligned} \label{eq:phitIm}
\end{align}
Hence, we have
\begin{align}
\phi(y_{-1})&=\Big(-\tfrac12-\tfrac{\sqrt3}2i,\,-\tfrac12+\tfrac{\sqrt3}2i,\,-i,\,i\Big) \notag \\
\phit(y_{-1})&=\Big(-\tfrac12+\tfrac{\sqrt3}2i,\,\tfrac12+\tfrac{\sqrt3}2i,\,
 \tfrac1{\sqrt2}-\tfrac1{\sqrt2}i,\,-\tfrac1{\sqrt2}-\tfrac1{\sqrt2}i\Big) \label{eq:phity-1} \\
\phi(x_0)&=\Big(-i,\,i,\,-\tfrac12-\tfrac{\sqrt3}2i,\,-\tfrac12+\tfrac{\sqrt3}2i\Big) \notag \\
\phit(x_0)&=\Big(-\tfrac1{\sqrt2}+\tfrac1{\sqrt2}i,\,\tfrac1{\sqrt2}+\tfrac1{\sqrt2}i,\,
 \tfrac12-\tfrac{\sqrt3}2i,\,-\tfrac12-\tfrac{\sqrt3}2i\Big). \label{eq:phitx0}
\end{align}
On $U_0$, for $0\le\theta\le\pi$,
\begin{align}
\phi(\tau_0(e^{i\theta}))&=\Big(e^{i(\theta+\frac\pi2)},\,e^{i(\theta-\frac\pi2)},\,
 -\tfrac12-\tfrac{\sqrt3}2i,\,-\tfrac12+\tfrac{\sqrt3}2i\Big) \notag \\
\phit(\tau_0(e^{i\theta}))&=\Big(e^{i(\frac\theta2+\frac\pi4)},\,e^{i(\frac\theta2-\frac\pi4)},\,
 \tfrac12-\tfrac{\sqrt3}2i,\,-\tfrac12-\tfrac{\sqrt3}2i\Big). \label{eq:phitU0}
\end{align}
In particular,
\begin{align}
\phi(y_0)&=\Big(i,\,-i,\,-\tfrac12-\tfrac{\sqrt3}2i,\,-\tfrac12+\tfrac{\sqrt3}2i\Big) \notag \\
\phit(y_0)&=\Big(\tfrac1{\sqrt2}+\tfrac1{\sqrt2}i,\,\tfrac1{\sqrt2}-\tfrac1{\sqrt2}i,\,
 \tfrac12-\tfrac{\sqrt3}2i,\,-\tfrac12-\tfrac{\sqrt3}2i\Big). \label{eq:phity0}
\end{align}
On $I_r$, for $t\in(0,1]$,
\begin{align}
\phi(\alpha(t))&=\Big(\tfrac t2+\tfrac{\sqrt{4-t^2}}2\,i,\,\tfrac t2-\tfrac{\sqrt{4-t^2}}2\,i,\,
\tfrac{-1-t}2-\tfrac{\sqrt{4-(1+t)^2}}2\,i,\,\tfrac{-1-t}2+\tfrac{\sqrt{4-(1+t)^2}}2\,i\Big) \notag \\[0.5ex]
\phit(\alpha(t))&=\begin{aligned}[t]\Big(&\sqrt{\tfrac12+\tfrac t4}+\sqrt{\tfrac12-\tfrac t4}\,i,\,
\sqrt{\tfrac12+\tfrac t4}-\sqrt{\tfrac12-\tfrac t4}\,i,\, \\
&\sqrt{\tfrac12-\tfrac{t+1}4}-\sqrt{\tfrac12+\tfrac{t+1}4}\,i,\,
-\sqrt{\tfrac12-\tfrac{t+1}4}-\sqrt{\tfrac12+\tfrac{t+1}4}\,i\Big).\end{aligned} \label{eq:phitIr}
\end{align}
In particular,
\begin{align}
\phi(r)&=\Big(\tfrac12+\tfrac{\sqrt3}2i,\,\tfrac12-\tfrac{\sqrt3}2i,-1,-1\Big) \notag \\
\phit(r)&=\Big(\tfrac{\sqrt3}2+\tfrac12i,\,\tfrac{\sqrt3}2-\tfrac12i,\,-i,\,-i\Big). \label{eq:phitr}
\end{align}
On $U_{-1}$, for $0\le\theta\le\pi$,
\begin{align}
\phi(\tau_{-1}(e^{i\theta}))&=\Big(-\tfrac12-\tfrac{\sqrt3}2i,\,-\tfrac12+\tfrac{\sqrt3}2i,\,
 e^{i(\theta-\frac\pi2)},\,e^{i(\theta+\frac\pi2)}\Big). \notag \\
\phit(\tau_{-1}(e^{i\theta}))&=\Big(-\tfrac12+\tfrac{\sqrt3}2i,\,\tfrac12+\tfrac{\sqrt3}2i,\,
 e^{i(\frac\theta2-\frac\pi4)},\,e^{i(\frac\theta2-\frac{3\pi}4)}\Big). \label{eq:phitU-1}
\end{align}
In particular,
\begin{align}
\phi(x_{-1})&=\Big(-\tfrac12-\tfrac{\sqrt3}2i,\,-\tfrac12+\tfrac{\sqrt3}2i,\,i,\,-i\Big) \notag \\
\phit(x_{-1})&=\Big(-\tfrac12+\tfrac{\sqrt3}2i,\,\tfrac12+\tfrac{\sqrt3}2i,\,
\tfrac1{\sqrt2}+\tfrac1{\sqrt2}i,\,\tfrac1{\sqrt2}-\tfrac1{\sqrt2}i\Big). \label{eq:phitx-1}
\end{align}
On $I_\ell$, for $t\in[-2,-1)$,
\begin{align}
\phi(\alpha(t))&=\Big(\tfrac t2-\tfrac{\sqrt{4-t^2}}2\,i,\,\tfrac t2+\tfrac{\sqrt{4-t^2}}2\,i,\,
\tfrac{-1-t}2+\tfrac{\sqrt{4-(1+t)^2}}2\,i,\,\tfrac{-1-t}2-\tfrac{\sqrt{4-(1+t)^2}}2\,i\Big) \notag \\[0.5ex]
\phit(\alpha(t))&=\begin{aligned}[t]\Big(&-\sqrt{\tfrac12+\tfrac t4}+\sqrt{\tfrac12-\tfrac t4}\,i,\,
\sqrt{\tfrac12+\tfrac t4}+\sqrt{\tfrac12-\tfrac t4}\,i,\, \\
&\sqrt{\tfrac12-\tfrac{t+1}4}+\sqrt{\tfrac12+\tfrac{t+1}4}\,i,\,
\sqrt{\tfrac12-\tfrac{t+1}4}-\sqrt{\tfrac12+\tfrac{t+1}4}\,i\Big).\end{aligned} \label{eq:phitIl}
\end{align}
In particular,
\begin{align}
\phi(s)&=\Big(-1,\,-1,\,\tfrac12+\tfrac{\sqrt3}2i,\,\tfrac12-\tfrac{\sqrt3}2i\Big) \notag \\
\phit(s)&=\Big(i,\,i,\,\tfrac{\sqrt3}2+\tfrac12i,\,\tfrac{\sqrt3}2-\tfrac12i\Big). \label{eq:phits}
\end{align}
\begin{prop}\label{prop:psigbe}
Let $b,b'\in\Bt$, $g,g'\in V$ and $\eps,\eps'\in E$.
Then
\begin{equation}\label{eq:psigbe}
\psi_g^\eps(b)=\psi_{g'}^{\eps'}(b')
\end{equation}
if and only if at least one of the following holds:
\renewcommand{\labelenumi}{(\alph{enumi})}
\begin{enumerate}

\item $b=b'$, $g=g'$ and $\eps=\eps'$

\item $b=b'=q$ and $\eps=\eps'$

\item $b=b'\in\partial_\ell\Bt$, $g^{-1}g'=t$ and $\eps=\eps'$

\item $b=b'\in\partial_r\Bt$, $g^{-1}g'=u$ and $\eps=\eps'$

\item $b=x_0$, $b'=y_0$, $(g,g')\in\{(e,t),(u,v)\}$ and $\eps\eps'=(-1,1,1,1)$

\item $b=x_0$, $b'=y_0$, $(g,g')\in\{(t,e),(v,u)\}$ and $\eps\eps'=(1,-1,1,1)$

\item $b=y_0$, $b'=x_0$, $(g,g')\in\{(t,e),(v,u)\}$ and $\eps\eps'=(-1,1,1,1)$

\item $b=y_0$, $b'=x_0$, $(g,g')\in\{(e,t),(u,v)\}$ and $\eps\eps'=(1,-1,1,1)$

\item $b=x_{-1}$, $b'=y_{-1}$, $(g,g')\in\{(e,u),(t,v)\}$ and $\eps\eps'=(1,1,-1,1)$

\item $b=x_{-1}$, $b'=y_{-1}$, $(g,g')\in\{(u,e),(v,t)\}$ and $\eps\eps'=(1,1,1,-1)$

\item $b=y_{-1}$, $b'=x_{-1}$, $(g,g')\in\{(u,e),(v,t)\}$ and $\eps\eps'=(1,1,-1,1)$

\item $b=y_{-1}$, $b'=x_{-1}$, $(g,g')\in\{(e,u),(t,v)\}$ and $\eps\eps'=(1,1,1,-1)$

\item $b=b'\in I_\ell$, $g^{-1}g'=v$ and $\eps\eps'=(-1,-1,1,1)$

\item $b=b'\in I_m$, $g^{-1}g'=v$ and $\eps\eps'=(-1,-1,-1,-1)$

\item $b=b'\in I_r$, $g^{-1}g'=v$ and $\eps\eps'=(1,1,-1,-1)$

\item $b=b'=s$, $g^{-1}g'=u$ and $\eps\eps'=(-1,-1,1,1)$

\item $b=b'=r$, $g^{-1}g'=t$ and $\eps\eps'=(1,1,-1,-1)$

\item $b=\tau_0(e^{i\theta})$, $b'=\tau_0(e^{i(\pi-\theta)})$ for some $0\le\theta\le\pi$,
      $(g,g')\in\{(e,u),(u,e)\}$ and $\eps\eps'=(-1,1,-1,-1)$

\item $b=\tau_0(e^{i\theta})$, $b'=\tau_0(e^{i(\pi-\theta)})$ for some $0\le\theta\le\pi$,
      $(g,g')\in\{(t,v),(v,t)\}$ and $\eps\eps'=(1,-1,-1,-1)$

\item $b=\tau_{-1}(e^{i\theta})$, $b'=\tau_{-1}(e^{i(\pi-\theta)})$ for some $0\le\theta\le\pi$,
      $(g,g')\in\{(e,t),(t,e)\}$ and $\eps\eps'=(-1,-1,1,-1)$

\item $b=\tau_{-1}(e^{i\theta})$, $b'=\tau_{-1}(e^{i(\pi-\theta)})$ for some $0\le\theta\le\pi$,
      $(g,g')\in\{(u,v),(v,u)\}$ and $\eps\eps'=(-1,-1,-1,1)$.

\end{enumerate}
\end{prop}
\begin{proof}
We will use the notation $\phit(b)=(z_1(b),z_2(b),z_3(b),z_4(b))$.
The equality~\eqref{eq:psigbe} holds if and only if either
\begin{equation}\label{eq:phitgbe}
\phit_g^\eps(b)=\phit_{g'}^{\eps'}(b')
\end{equation}
or
\begin{equation}\label{eq:phitgbebar}
\phit_g^\eps(b)=\overline{\phit_{g'}^{\eps'}(b')}
\end{equation}
holds.

Suppose~\eqref{eq:phitgbe} holds.
Then squaring yields $\phi_g(b)=\phi_{g'}(b')$,
and by Lemma~\ref{lem:phiB},
at least one of conditions (i)--(vi) listed there holds.
Clearly (i)$\Rightarrow$(a) and (ii)$\Rightarrow$(b).
We have $z_1(b)=z_2(b)$ for all $b\in\partial_\ell\Bt$, by continuity and since it holds
at $b=q$, and thus (iii)$\Rightarrow$(c).
Similarly, since $z_3(b)=z_4(b)$ for all $b\in\partial_r\Bt$, we get (iv)$\Rightarrow$(d).
Conversely, it is clear that any of (a)--(d) imply~\eqref{eq:phitgbe}.
Suppose~(v) holds.
Since $b,b'\in\Bt$, we need $\zeta=\pm1$, i.e. either $b=x_0$ and $b'=y_0$ or $b=y_0$ and $b'=x_0$.
Examining~\eqref{eq:phitx0} and~\eqref{eq:phity0}, we obtain that one of (e)--(h) holds,
and conversely any of (e)--(h) implies~\eqref{eq:phitgbe}.
Similarly, (vi) and~\eqref{eq:phitgbe} together are equivalent to (i)--(l).

Now suppose~\eqref{eq:phitgbebar} holds.
Squaring and using~\eqref{eq:phibeta}, we get
\begin{equation}\label{eq:phigbeta}
\phi_g(b)=\overline{\phi_{g'}(b')}=\phi_{g'v}(\beta(b')).
\end{equation}
From Lemma~\ref{lem:phiB}, we get that at least one of the following holds:
\renewcommand{\labelenumi}{(\roman{enumi}$'$)}
\begin{enumerate}

\item $b=\beta(b')$ and $g^{-1}g'=v$

\item $b=\beta(b')\in\partial_\ell B\cap\partial_r B$

\item $b=\beta(b')\in\partial_\ell B$ and $g^{-1}g'=u$

\item $b=\beta(b')\in\partial_r B$ and $g^{-1}g'=t$

\item $b=\tau_0(\zeta)$ and $b'=\beta(\tau_0(-\zeta))=\tau_0(-\overline\zeta)$ for
      some $\zeta\in\Tcirc$, and $g^{-1}g'=u$

\item $b=\tau_{-1}(\zeta)$ and $b'=\beta(\tau_{-1}(-\zeta))=\tau_{-1}(-\overline\zeta)$ for
      some $\zeta\in\Tcirc$, and $g^{-1}g'=t$.

\end{enumerate}
Since $b,b'\in\Bt$, in cases (i$'$)--(iv$'$), we have
$b=\beta(b')$, so $b\in\Bt\cap\beta(\Bt)=I_\ell\cup I_m\cup I_r$.
But $\beta$ restricts to the identity map on $\Bt\cap\beta(\Bt)$, and we conclude $b=b'$.
In particular condition~(ii$'$) cannot hold.

Suppose~(i$'$) holds.
Then~\eqref{eq:phitgbebar} becomes
\begin{equation}\label{eq:phitgv}
\phit_g^\eps(b)=\overline{\phit_{gv}^{\eps'}(b)}.
\end{equation}
If $b\in I_\ell$, then examining~\eqref{eq:phitx-1}, \eqref{eq:phitIl} and~\eqref{eq:phits},
we conclude that~\eqref{eq:phitgv} holds if and only if $\eps\eps'=(-1,-1,1,1)$;
this corresponds to condition~(m).
If $b\in I_m$, then examining~\eqref{eq:phitIm}, \eqref{eq:phity-1} and~\eqref{eq:phitx0},
we find that~\eqref{eq:phitgv} holds if and only if $\eps\eps'=(-1,-1,-1,-1)$;
this corresponds to condition~(n).
If $b\in I_r$, then examining~\eqref{eq:phity0}, \eqref{eq:phitIr} and~\eqref{eq:phitr},
we find that~\eqref{eq:phitgv} holds if and only if $\eps\eps'=(1,1,-1,-1)$;
this corresponds to condition~(o).

Suppose~(iii$'$) holds.
Then $b=b'=s$.
Using~\eqref{eq:phits}, we conclude that~\eqref{eq:phitgv} holds if and only if $\eps\eps'=(-1,-1,1,1)$;
this corresponds to condition~(p).

Suppose~(iv$'$) holds.
Then $b=b'=r$.
Using~\eqref{eq:phitr}, we conclude that~\eqref{eq:phitgv} holds if and only if $\eps\eps'=(1,1,-1,-1)$;
this corresponds to condition~(q).

Suppose~(v$'$) holds, with $b=\tau_0(e^{i\theta})$, $0\le\theta\le\pi$;
we have $b'=\tau_0(e^{i(\pi-\theta)})$.
Then from~\eqref{eq:phitU0},
\begin{align*}
\phit(b)&=\Big(e^{i(\frac\theta2+\frac\pi4)},\,e^{i(\frac\theta2-\frac\pi4)},\,
 \tfrac12-\tfrac{\sqrt3}2i,\,-\tfrac12-\tfrac{\sqrt3}2i\Big) \\
\overline{\phit(b')}&=\Big(e^{i(\frac\theta2-\frac{3\pi}4)},\,e^{i(\frac\theta2-\frac\pi4)},\,
 \tfrac12+\tfrac{\sqrt3}2i,\,-\tfrac12+\tfrac{\sqrt3}2i\Big).
\end{align*}
Thus~\eqref{eq:phitgv} holds if and only if~(r) or~(s) holds.

Similarly, if~(vi$'$) holds, then~\eqref{eq:phitgv} holds if and only if~(t) or~(u) holds.
\end{proof}

A consequence of Proposition~\ref{prop:psigbe} is that each of the sixty--four maps
$\psi_g^\eps:\Bt\to\Ec/\Ints_2$ is injective, and is thus a homeomorphism onto its image.
We will use the notation $\Bt_g^\eps$ for the image of $\psi_g^\eps$, identified with $\Bt$
via $\psi_g^\eps$.
The space $\Ec/\Ints_2$ is homeomorphic to the identification space obtained
by gluing together these sixty--four pieces $(\Bt_g^\eps)_{g\in V,\,\eps\in E}$
according to how they overlap in $\Ec/\Ints_2$, namely as explicated in
Proposition~\ref{prop:psigbe}.

Given $\eps\in E$, we glue together $\Bt_e^\eps$, $\Bt_t^\eps$, $\Bt_u^\eps$ and $\Bt_v^\eps$
via (b)--(d) of Proposition~\ref{prop:psigbe}
to obtain the solid square piece $\Bt^\eps$ depicted in Figure~\ref{fig:Bteps}.
(Note, however, that the labels on the points of the boundary of $\Bt^\eps$ are those
retained from the picture of $\Bt$, and do not indicate that identifications of these point
are made;
thus, for example, the two points labeled ``$s$'' in Figure~\ref{fig:Bteps} are not identified
with each other.)
\begin{figure}[bp]
\caption{The square $\Bt^\eps$ assembled from $(\Bt_g^\eps)_{g\in V}$.}\label{fig:Bteps}
\setlength{\unitlength}{0.97in}
\begin{picture}(6.0,6.0)(-0.5,-0.5)
\put(2.5,2.5){\circle*{0.05}}
\put(2.6,2.47){$q$}
\dashline[-20]{0.06}(0,0)(5,5)
\dashline[-20]{0.06}(0,5)(5,0)
\put(2.45,1.2){$\Bt_e^\eps$}
\put(1.1,0.9){$\partial_\ell\Bt_e^\eps$}
\put(3.55,0.9){$\partial_r\Bt_e^\eps$}
\put(0,0){\line(1,0){5}}
\multiput(0,0)(1,0){6}{\circle*{0.05}}
\put(-0.15,-0.15){$s$}
\put(0.9,-0.2){$x_{-1}$}
\put(1.9,-0.2){$y_{-1}$}
\put(2.93,-0.2){$x_0$}
\put(3.93,-0.2){$y_0$}
\put(5.1,-0.15){$r$}
\put(0.5,0.05){$I_\ell$}
\put(1.45,0.05){$U_{-1}$}
\put(2.45,0.05){$I_m$}
\put(3.45,0.05){$U_0$}
\put(4.4,0.05){$I_r$}
\put(3.6,2.4){$\Bt_u^\eps$}
\put(3.9,1.2){$\partial_r\Bt_u^\eps$}
\put(3.9,3.7){$\partial_\ell\Bt_u^\eps$}
\put(5,0){\line(0,1){5}}
\multiput(5,1)(0,1){5}{\circle*{0.05}}
\put(5.1,0.95){$y_0$}
\put(5.1,1.95){$x_0$}
\put(5.1,2.95){$y_{-1}$}
\put(5.1,3.95){$x_{-1}$}
\put(5.1,5.05){$s$}
\put(4.8,0.5){$I_r$}
\put(4.77,1.45){$U_0$}
\put(4.8,2.45){$I_m$}
\put(4.72,3.45){$U_{-1}$}
\put(4.8,4.4){$I_\ell$}
\put(1.3,2.4){$\Bt_t^\eps$}
\put(0.8,1.2){$\partial_\ell\Bt_t^\eps$}
\put(0.8,3.7){$\partial_r\Bt_t^\eps$}
\put(0,0){\line(0,1){5}}
\multiput(0,1)(0,1){5}{\circle*{0.05}}
\put(-0.3,0.95){$x_{-1}$}
\put(-0.3,1.95){$y_{-1}$}
\put(-0.23,2.97){$x_0$}
\put(-0.23,3.97){$y_0$}
\put(-0.15,5.1){$r$}
\put(0.05,0.5){$I_\ell$}
\put(0.05,1.45){$U_{-1}$}
\put(0.05,2.45){$I_m$}
\put(0.05,3.45){$U_0$}
\put(0.05,4.4){$I_r$}
\put(2.45,3.7){$\Bt_e^\eps$}
\put(1.1,4){$\partial_r\Bt_e^\eps$}
\put(3.55,4){$\partial_\ell\Bt_e^\eps$}
\put(0,5){\line(1,0){5}}
\multiput(1,5)(1,0){4}{\circle*{0.05}}
\put(0.95,5.13){$y_0$}
\put(1.95,5.13){$x_0$}
\put(2.95,5.13){$y_{-1}$}
\put(3.95,5.13){$x_{-1}$}
\put(0.5,4.8){$I_r$}
\put(1.45,4.8){$U_0$}
\put(2.45,4.8){$I_m$}
\put(3.45,4.8){$U_{-1}$}
\put(4.4,4.8){$I_\ell$}
\end{picture}
\end{figure}
The remaining parts (e)--(u) of Proposition~\ref{prop:psigbe}
are instructions for gluing the sixteen squares $(\Bt^\eps)_{\eps\in E}$ along certain of the
edges and vertices, in order to obtain $\Ec/\Ints_2$.
In order to describe this space, we relabel the edges and vertices of $\Bt^\eps$
as shown in Figure~\ref{fig:Btepslab}.
\begin{figure}[bp]
\caption{$\Bt^\eps$ relabeled, with edges oriented.}\label{fig:Btepslab}
\setlength{\unitlength}{0.97in}
\begin{picture}(6.0,6.0)(-0.5,-0.5)
\put(0,0){\line(1,0){5}}
\multiput(0,0)(1,0){6}{\circle*{0.05}}
\multiput(0.575,0)(1,0){5}{\line(-2,1){0.15}}
\multiput(0.575,0)(1,0){5}{\line(-2,-1){0.15}}
\put(-0.3,-0.2){$a(\eps)$}
\multiputlist(1,-0.2)(1,0){$b(\eps)$,$c(\eps)$,$d(\eps)$,$e(\eps)$}
\put(5,-0.2){$f(\eps)$}
\multiputlist(0.5,0.22)(1,0){$A(\eps)$,$B(\eps)$,$C(\eps)$,$D(\eps)$,$E(\eps)$}
\put(5,0){\line(0,1){5}}
\multiput(5,1)(0,1){5}{\circle*{0.05}}
\multiput(5,0.575)(0,1){5}{\line(1,-2){0.075}}
\multiput(5,0.575)(0,1){5}{\line(-1,-2){0.075}}
\multiputlist(5.25,1)(0,1){$g(\eps)$,$h(\eps)$,$i(\eps)$,$j(\eps)$}
\put(5.1,5.1){$k(\eps)$}
\multiputlist(4.7,0.5)(0,1){$F(\eps)$,$G(\eps)$,$H(\eps)$,$I(\eps)$,$J(\eps)$}
\put(5,5){\line(-1,0){5}}
\multiput(4,5)(-1,0){5}{\circle*{0.05}}
\multiput(4.425,5)(-1,0){5}{\line(2,-1){0.15}}
\multiput(4.425,5)(-1,0){5}{\line(2,1){0.15}}
\multiputlist(4,5.17)(-1,0){$\ell(\eps)$,$m(\eps)$,$n(\eps)$,$o(\eps)$}
\put(-0.28,5.13){$p(\eps)$}
\multiputlist(4.5,4.75)(-1,0){$K(\eps)$,$L(\eps)$,$M(\eps)$,$N(\eps)$,$O(\eps)$}
\put(0,5){\line(0,-1){5}}
\multiput(0,4)(0,-1){4}{\circle*{0.05}}
\multiput(0,4.425)(0,-1){5}{\line(1,2){0.075}}
\multiput(0,4.425)(0,-1){5}{\line(-1,2){0.075}}
\multiputlist(-0.2,4)(0,-1){$q(\eps)$,$r(\eps)$,$s(\eps)$,$t(\eps)$}
\multiputlist(0.27,4.5)(0,-1){$P(\eps)$,$Q(\eps)$,$R(\eps)$,$S(\eps)$,$T(\eps)$}
\end{picture}
\end{figure}
The vertices are lowercase, the edges uppercase, and we have oriented the edges as shown.
We will now describe the identifications of vertices and edges that occur;
all identifications of edges are orientation preserving.
From~(m), identify
\begin{alignat*}{2}
\mbox{edges }&&A(\eps)&\mbox{ with }K((-1,-1,1,1)\eps) \\
&&J(\eps)&\mbox{ with }T((-1,-1,1,1)\eps) \\
\mbox{vertices }&&a(\eps)&\mbox{ with }k((-1,-1,1,1)\eps) \\
&&b(\eps)&\mbox{ with }\ell((-1,-1,1,1)\eps) \\
&&j(\eps)&\mbox{ with }t((-1,-1,1,1)\eps).
\end{alignat*}
From~(n), identify
\begin{alignat*}{2}
\mbox{edges }&&C(\eps)&\mbox{ with }M((-1,-1,-1,-1)\eps) \\
&&H(\eps)&\mbox{ with }R((-1,-1,-1,-1)\eps) \\
\mbox{vertices }&&c(\eps)&\mbox{ with }m((-1,-1,-1,-1)\eps) \\
&&d(\eps)&\mbox{ with }n((-1,-1,-1,-1)\eps) \\
&&h(\eps)&\mbox{ with }r((-1,-1,-1,-1)\eps) \\
&&i(\eps)&\mbox{ with }s((-1,-1,-1,-1)\eps).
\end{alignat*}
From~(o), identify
\begin{alignat*}{2}
\mbox{edges }&&E(\eps)&\mbox{ with }O((1,1,-1,-1)\eps) \\
&&F(\eps)&\mbox{ with }P((1,1,-1,-1)\eps) \\
\mbox{vertices }&&e(\eps)&\mbox{ with }o((1,1,-1,-1)\eps) \\
&&f(\eps)&\mbox{ with }p((1,1,-1,-1)\eps) \\
&&g(\eps)&\mbox{ with }q((1,1,-1,-1)\eps).
\end{alignat*}
Part~(p) identifies vertices $a(\eps)$ with $k((-1,-1,1,1)\eps)$, which has already been done.
Part~(q) identifies vertices $f(\eps)$ with $p((1,1,-1,-1)\eps)$, which has already been done.
From~(r), identify
\begin{alignat*}{2}
\mbox{edges }&&D(\eps)&\mbox{ with }G((-1,1,-1,-1)\eps) \\
\mbox{vertices }&&d(\eps)&\mbox{ with }g((-1,1,-1,-1)\eps) \\
&&e(\eps)&\mbox{ with }h((-1,1,-1,-1)\eps).
\end{alignat*}
From~(s), identify
\begin{alignat*}{2}
\mbox{edges }&&N(\eps)&\mbox{ with }Q((1,-1,-1,-1)\eps) \\
\mbox{vertices }&&n(\eps)&\mbox{ with }q((1,-1,-1,-1)\eps) \\
&&o(\eps)&\mbox{ with }r((1,-1,-1,-1)\eps).
\end{alignat*}
From~(t), identify
\begin{alignat*}{2}
\mbox{edges }&&B(\eps)&\mbox{ with }S((-1,-1,1,-1)\eps) \\
\mbox{vertices }&&b(\eps)&\mbox{ with }s((-1,-1,1,-1)\eps) \\
&&c(\eps)&\mbox{ with }t((-1,-1,1,-1)\eps).
\end{alignat*}
From~(u), identify
\begin{alignat*}{2}
\mbox{edges }&&I(\eps)&\mbox{ with }L((-1,-1,-1,1)\eps) \\
\mbox{vertices }&&i(\eps)&\mbox{ with }\ell((-1,-1,-1,1)\eps) \\
&&j(\eps)&\mbox{ with }m((-1,-1,-1,1)\eps).
\end{alignat*}
The remaining parts of Proposition~\ref{prop:psigbe} make identifications that have already been done.
Namely, part~(e) identifies vertices $d(\eps)$ with $q((-1,1,1,1)\eps)$
and $h(\eps)$ with $o((-1,1,1,1)\eps)$;
part~(f) identifies vertices $r(\eps)$ with $e((1,-1,1,1)\eps)$
and $n(\eps)$ with $g((1,-1,1,1)\eps)$;
parts~(g) and~(h) simply repeat parts~(e) and~(f);
part~(i) identifies vertices $b(\eps)$ with $i((1,1,-1,1)\eps)$
and $t(\eps)$ with $m((1,1,-1,1)\eps)$;
part~(j) identifies vertices $j(\eps)$ with $c((1,1,1,-1)\eps)$
and $\ell(\eps)$ with $s((1,1,1,-1)\eps)$;
parts~(k) and~(l) repeat parts~(i) and~(j).

Therefore, the identification space, which is homeomorphic to $\Ec/\Ints_2$, and thus to $\Gc^\Reals_{5,2}$,
is a 2--complex having sixteen 2--cells, namely the squares $(\Bt^\eps)_{\eps\in E}$, one hundred and sixty
edges,
\[
\bigcup_{\eps\in E}\{\Abar(\eps),\Bbar(\eps),\Cbar(\eps),\Dbar(\eps),\Ebar(\eps),
\Fbar(\eps),\Hbar(\eps),\Ibar(\eps),\Jbar(\eps),\Nbar(\eps)\},
\]
and ninety--six vertices,
\[
\bigcup_{\eps\in E}\{\abar(\eps),\bbar(\eps),\cbar(\eps),\dbar(\eps),\ebar(\eps),\fbar(\eps)\},
\]
where
\begin{align*}
\Abar(\eps)&=\{A(\eps),K((-1,-1,1,1)\eps)\} \\
\Bbar(\eps)&=\{B(\eps),S((-1,-1,1,-1)\eps)\} \\
\Cbar(\eps)&=\{C(\eps),M((-1,-1,-1,-1)\eps)\} \\
\Dbar(\eps)&=\{D(\eps),G((-1,1,-1,-1)\eps)\} \\
\Ebar(\eps)&=\{E(\eps),O((1,1,-1,-1)\eps)\} \\
\Fbar(\eps)&=\{F(\eps),P((1,1,-1,-1)\eps)\} \\
\Hbar(\eps)&=\{H(\eps),R((-1,-1,-1,-1)\eps)\} \\
\Ibar(\eps)&=\{I(\eps),L((-1,-1,-1,1)\eps)\} \\
\Jbar(\eps)&=\{J(\eps),T((-1,-1,1,1)\eps)\} \\
\Nbar(\eps)&=\{N(\eps),Q((1,-1,-1,-1)\eps)\}
\end{align*}
and where
\begin{align*}
\abar(\eps)&=\{a(\eps),k((-1,-1,1,1)\eps)\} \\
\bbar(\eps)&=\{b(\eps),i((1,1,-1,1)\eps),\ell((-1,-1,1,1)\eps),s((-1,-1,1,-1)\eps)\} \\
\cbar(\eps)&=\{c(\eps),j((1,1,1,-1)\eps),m((-1,-1,-1,-1)\eps),t((-1,-1,1,-1)\eps)\} \\
\dbar(\eps)&=\{d(\eps),g((-1,1,-1,-1)\eps),n((-1,-1,-1,-1)\eps),q((-1,1,1,1)\eps)\} \\
\ebar(\eps)&=\{e(\eps),h((-1,1,-1,-1)\eps),o((1,1,-1,-1)\eps),r((1,-1,1,1)\eps)\} \\
\fbar(\eps)&=\{f(\eps),p((1,1,-1,-1)\eps)\}.
\end{align*}

From Theorem~\ref{thm:mfld}, $\Gc^\Reals_{5,2}$ is a closed, orientable surface.
From the above description, the Euler characteristic of $\Gc^\Reals_{5,2}$ is $16-160+96=-48$.

\begin{thm}\label{thm:G52}
The space $\Gc^\Reals_{5,2}$ is homeomorphic to the closed, orientable surface
of genus $25$.
\end{thm}

\section{Connectedness}\label{sec:connected}

Since $\Gc^\Eb_{k,n}$ and $\Fc^\Eb_{k,n}$ are real algebraic sets, each has only
finitely many connected components, by the classical result of Whitney~\cite{W}.
We already saw, in~\S\ref{sec:1red}, that $\Gc^\Cpx_{n+1,1}$ and thus $\Gc^\Cpx_{n+1,n}$
are connected, while $\Gc^\Reals_{n+1,1}$ and $\Gc^\Reals_{n+1,n}$ are disconnected, for $n\ge1$.

Recall from Corollary~\ref{cor:stf}
we have the locally trivial fiber bundle $\pi^\Eb_{k,n}:\Fc^\Eb_{k,n}\to\Gc^\Eb_{k,n}$
and $\Gc^\Eb_{k,n}$ is thereby identified with the orbit space $\Fc^\Eb_{k,n}/\Oc^\Eb_n$.
Moreover, from Corollary~\ref{cor:I-P}, we have the homeomorphism $\gamma_{k,n}:\Gc^\Eb_{k,n}\to\Gc^\Eb_{k,k-n}$.

In the complex case, since the fibers $\Oc_n^\Cpx$ are connected, the following result is obvious.
\begin{prop}
Let $k,n\in\Nats$ with $k>n$.
Let $C$ be a connected component of $\Gc^\Cpx_{k,n}$.
Then $(\pi_{k,n}^\Cpx)^{-1}(C)$ is connected.
Thus, $\Fc^\Cpx_{k,n}$ and $\Gc^\Cpx_{k,n}$ have the same number of connected components;
in particular, $\Fc^\Cpx_{k,n}$ is connected if and only if $\Gc^\Cpx_{k,n}$ is connected.
Furthermore, $\Fc^\Cpx_{k,n}$ is connected if and only if $\Fc^\Cpx_{k,k-n}$ is connected.
\end{prop}

The real case, however, is somewhat more interesting.
\begin{prop}\label{prop:connected}
Let $k,n\in\Nats$ with $k>n$.
Let $C$ be a connected component of $\Gc^\Reals_{k,n}$.
Then $(\pi^\Reals_{k,n})^{-1}(C)$ is either connected or has exactly two connected components.
Furthermore, $(\pi^\Reals_{k,n})^{-1}(C)$ and $(\pi^\Reals_{k,k-n})^{-1}(\gamma_{k,n}(C))$ have the same number
of connected components.
\end{prop}
\begin{proof}
Let $\Fc^\Reals_{k,n}/\SOc_n$ denote the orbit space of the restriction of the action of $\Oc^\Reals_n$
on $\Fc^\Reals_{k,n}$ to the special orthogonal group $\SOc_n\subset\Oc^\Reals_n$.
Since $\pi_{k,n}^\Reals$ is a locally trivial
fiber bundle, we see that the quotient $q:\Fc^\Reals_{k,n}/\SOc_n\to\Gc^\Reals_{k,n}$
is a two--fold covering projection, and that $(\pi^\Reals_{k,n})^{-1}(C)$ and $q^{-1}(C)$ have the same
number of connected components; in particular this number is either one or two.

Since $\Gc^\Reals_{k,n}$ is a real algebraic variety, by Whitney's results~\cite{W} (see also \S\ref{sec:mfld})
it is locally path connected and thus $C$ is path connected.
If $(\pi^\Reals_{k,n})^{-1}(C)$ is connected, then there is a closed path $\tau:[0,1]\to C$ with the property that
if $\taut:[0,1]\to\Fc^\Reals_{k,n}$ is a lifting, then letting $U\in\Oc_n$ be such that $\taut(1)=U\taut(0)$,
we have $\det(U)=-1$.
Let $\taut':[0,1]\to\Fc^\Reals_{k,k-n}$ be a lifting of $\gamma_{k,n}\circ\tau:[0,1]\to\Gc^\Reals_{k,k-n}$
and let $U'\in\Oc_{k-n}$ be such that $\taut'(1)=U'\taut'(0)$.
We will show $\det(U')=-1$, which will imply $(\pi^\Reals_{k,k-n})^{-1}(C)$ is connected and will
thus finish the proof.
Let $P(t)=\frac nk\tau(t)$.
Then, {\em cf}~\eqref{eq:FWU},
\[
V(t)=\sqrt{\frac nk}W_{n,k}^*\taut(t)\in M_k(\Reals),\quad(t\in[0,1])
\]
is a continuous path of partial isometries satisfying
\[
V(t)^*V(t)=P(t),\qquad V(t)V(t)^*=\diag(\underset{n}{\underbrace{1,\ldots,1}},0,\ldots,0),\qquad(t\in[0,1])
\]
and $V(1)V(0)^*=\left(\begin{smallmatrix}U&0\\ 0&0\end{smallmatrix}\right)$.
Similarly,
\[
V'(t)=\sqrt{\frac{k-n}n}\left(\begin{matrix}0_{n,k-n}\\I_{k-n}\end{matrix}\right)\taut'(t)
\in M_k(\Reals),\quad(t\in[0,1])
\]
is a continuous path of partial isometries satisfying
\[
V'(t)^*V'(t)=I_k-P(t),\qquad V'(t)V'(t)^*=\diag(0,\ldots,0,\underset{k-n}{\underbrace{1,\ldots,1}}),
\qquad(t\in[0,1])
\]
and $V'(1)V'(0)^*=\left(\begin{smallmatrix}0&0\\ 0&U'\end{smallmatrix}\right)$.
Therefore,
$(V(t)+V'(t))(V(0)+V'(0))^*$, $0\le t\le1$, is a continuous path in $\Oc^\Reals_k$ starting at $I_k$ and ending
at $\left(\begin{smallmatrix}U&0\\ 0&U'\end{smallmatrix}\right)$.
This implies $\det(U)=\det(U')$.
\end{proof}

\begin{cor}\label{cor:connected}
Let $k,n\in\Nats$, $k>n$.
Then $\Fc^\Reals_{k,n}$ is connected if and only if $\Fc^\Reals_{k,k-n}$ is connected.
\end{cor}

\begin{thm}\label{thm:Fk2}
Let $k\in\Nats$, $k\ge4$.
Then $\Fc^\Reals_{k,2}$ is connected.
\end{thm}
\begin{proof}
By Corollary~\ref{cor:Fct}, $\Fc^\Reals_{k,2}$ is homeomorphic to
\[
\Fct_{k,2}=\{(z_1,\ldots,z_k)\in\Tcirc^k\mid\sum_{j=1}^kz_j^2=0\}.
\]
Let
\[
\Cc=\{(w_1,\ldots,w_k)\in\Tcirc^k\mid\sum_{j=1}^kw_j=0\}
\]
and let $p:\Fct_{k,2}\to\Cc$ be $p((z_1,\ldots,z_k))=(z_1^2,\ldots,z_k^2)$.
Then $p$ is a $2^k$--fold covering projection.

Thus $\Cc$ is the space of all chains in $\Reals^2$, starting and ending at $0$ and
having links of uniform length $1$.
It is known (see, for example,~\cite[Thm. 3.1]{LW}) that the space $\Cc$ is path connected.
Thus, given any $c\in\Cc$, there is a path from $c$ to some chosen element $s\in\Cc$, that is
said to be in {\em standard form} and is described below.
Of course, given a path $\gamma:[0,1]\to\Cc$ and given $a\in p^{-1}(\gamma(0))$, there is
a (unique) path $\gamma':[0,1]\to\Fct_{k,2}$ such that $p\circ\gamma'=\gamma$ and $\gamma'(0)=a$.
Therefore, in order to show connectedness of $\Fct_{k,2}$, it will suffice to specify a particular
element $b\in p^{-1}(s)$ and to exhibit a path in $\Fct_{k,2}$ from each element of $p^{-1}(s)$ to $b$.

\noindent{\bf Case I:} $k=4$.
We take as standard element $s=(1,-1,1,-1)\in\Cc$.
Then
\[
p^{-1}(s)=\{(\eps_1,\eps_2i,\eps_3,\eps_4i)\mid\eps_j\in\{\pm1\}\},
\]
and we select $b=(1,i,1,i)$.
It is not difficult to construct paths in $\Fct_{4,2}$ from all elements of $p^{-1}(s)$ to $b$.
This is somewhat tedious and is left to the reader.

Alternatively, it is proved in~\S\ref{sec:G42} that $\Gc^\Reals_{4,2}$ is connected;
the points $b$ and $a=(1,-i,1,-i)$ in $\Fct_{4,2}$ belong to the same $\Oc^\Reals_2$--orbit but differ
by a matrix in $\Oc^\Reals_2$ of determinant $-1$;
therefore, by the technique of the proof of Proposition~\ref{prop:connected},
in order to show that $\Fct_{4,2}$ is connected, it will suffice to find a path
in $\Fct_{4,2}$ from $b$ to $a$.
Starting at $b$ the path 
\[
(e^{i\theta},ie^{i\theta},1,i),\qquad0\le\theta\le\pi
\]
takes us to $(-1,-i,1,i)$;
then the path
\[
(e^{i\theta},-i,1,ie^{i\theta}),\qquad0\le\theta\le\pi,
\]
takes us to $(1,-i,1,-i)=a$.

\noindent{\bf Case II:} $k$ even, $k\ge6$.
Take as standard element $s=(1,-1,1,-1,\ldots,1,-1)\in\Cc$, so
\[
p^{-1}(s)=\{(\eps_1,\eps_2i,\eps_3,\eps_4i,\ldots,\eps_{k-1},\eps_ki)\mid\eps_j\in\{\pm1\}\}.
\]
Take $b=(1,i,\ldots,1,i)$ and, given
$a=(\eps_1,\eps_2i,\ldots,\eps_{k-1},\eps_ki)\in p^{-1}(s)$, construct a path in $\Fct_{k,2}$ from
$a$ to $b$ as follows.
We have $(\eps_1,\eps_2i,\eps_3,\eps_4i)\in\Fct_{4,2}$, and by Case~I, $\Fct_{4,2}$ is connected,
so there is a path in $\Fct_{4,2}$
from $(\eps_1,\eps_2i,\eps_3\eps_4i)$ to $(1,i,1,i)$.
Keeping the remaining $k-4$ elements constant, this yields a path in $\Fct_{k,2}$ from $a$ to
\begin{equation}\label{eq:a2}
(1,i,1,i,\eps_5,\eps_6i,\ldots,\eps_{k-1},\eps_{k-1}i).
\end{equation}
Now taking a path in $\Fct_{4,2}$ from $(1,i,\eps_5,\eps_6i)$ to $(1,i,1,i)$ yields
a path in $\Fct_{k,2}$ from the point in~\eqref{eq:a2} to
\[
(1,i,1,i,1,i,\eps_7,\eps_8i,\ldots,\eps_{k-1},\eps_{k-1}i).
\]
Continuing in this manner, we construct a path in $\Fct_{k,2}$ from $a$ to $b$.

\noindent{\bf Case III:} $k=5$.
Take as standard element
\[
s=(-\tfrac12+\tfrac{\sqrt3}2i,-\tfrac12-\tfrac{\sqrt3}2i,1,1,-1)\in\Cc.
\]
Then
\[
p^{-1}(s)=\{(\eps_1(\tfrac12+\tfrac{\sqrt3}2i),\eps_2(\tfrac12-\tfrac{\sqrt3}2i),\eps_3,\eps_4,\eps_5i)
\mid\eps_j\in\{\pm1\}\}
\]
and we select
\[
b=(\tfrac12+\tfrac{\sqrt3}2i,\tfrac12-\tfrac{\sqrt3}2i,1,1,i)=(e^{\pi i/3},e^{-\pi i/3},1,1,i).
\]
It is again routine, though tedious, to construct paths in $\Fct_{5,2}$ from all the elements of $p^{-1}(s)$
to $b$.

Alternatively, arguing as in Case~I above, using the result from~\S\ref{sec:G52} that $\Gc^\Reals_{5,2}$
is connected, it will suffice to construct a path in $\Fct_{5,2}$ from $b$ to
\[
a=(\tfrac12-\tfrac{\sqrt3}2i,\tfrac12+\tfrac{\sqrt3}2i,1,1,-i)=(e^{-\pi i/3},e^{\pi i/3},1,1,-i).
\]
The path
\[
(e^{\pi i/3},e^{-\pi i/3},e^{i\theta},1,ie^{i\theta}),\qquad 0\le\theta\le\pi/3
\]
takes us from $b$ to $(e^{\pi i/3},e^{-\pi i/3},e^{\pi i/3},1,e^{5\pi i/6})$;
the path
\[
(e^{i(\theta+\pi/3)},e^{-\pi i/3},e^{\pi i/3},1,e^{i(\theta+5\pi/6)}),\qquad0\le\theta\le4\pi/3
\]
takes us to $(e^{-\pi i/3},e^{-\pi i/3},e^{\pi i/3},1,e^{\pi i/6})$;
the path
\[
(e^{-\pi i/3},e^{i(\theta-\pi/3)},e^{\pi i/3},1,e^{i(\theta+\pi/6)}),\qquad0\le\theta\le\pi/6
\]
takes us to $(e^{-\pi i/3},e^{-\pi i/6},e^{\pi i/3},1,e^{\pi i/3})$;
the path
\[
(e^{-\pi i/3},e^{i(\theta-\pi/6)},e^{i(\theta+\pi/3)},1,e^{\pi i/3}),\qquad0\le\theta\le\pi/2
\]
takes us to $(e^{-\pi i/3},e^{\pi i/3},e^{5\pi i/6},1,e^{\pi i/3})$;
the path
\[
(e^{-\pi i/3},e^{\pi i/3},e^{i(\theta+5\pi/6)},1,e^{i(\theta+\pi/3)}),\qquad0\le\theta\le7\pi/6
\]
takes us to $a$.

\noindent{\bf Case IV:} $k$ odd, $k\ge5$.
Take as standard element
\[
s=(-\tfrac12+\tfrac{\sqrt3}2i,-\tfrac12-\tfrac{\sqrt3}2i,1,1,-1,1,-1,\ldots,1,-1)\in\Cc.
\]
and select
\[
b=(\tfrac12+\tfrac{\sqrt3}2i,\tfrac12-\tfrac{\sqrt3}2i,1,1,i)=(e^{\pi i/3},e^{-\pi i/3},1,1,i,1,i,\ldots,1,i)
\in p^{-1}(s).
\]
We must find a path in $\Fct_{k,2}$ from an arbitrary element
\[
c=(\eps_1(\tfrac12+\tfrac{\sqrt3}2i),\eps_2(\tfrac12-\tfrac{\sqrt3}2i),\eps_3,\eps_4,\eps_5i,
\ldots,\eps_{k-1},\eps_ki),
\qquad\eps_j\in\{\pm1\}
\]
of $p^{-1}(s)$ to $b$.
Similarly to in Case~II above, using paths in $\Fct_{5,2}$, which by Case~III
we know to be connected,
we construct a path from $c$ to $b$
passing through the points
\begin{align*}
&(\tfrac12+\tfrac{\sqrt3}2i,\tfrac12-\tfrac{\sqrt3}2i,1,1,i,\eps_6,\eps_7i,\eps_8,\eps_9i,\ldots,
\ldots,\eps_{k-1},\eps_ki) \\
&(\tfrac12+\tfrac{\sqrt3}2i,\tfrac12-\tfrac{\sqrt3}2i,1,1,i,1,i,\eps_8,\eps_9i,\ldots,
\ldots,\eps_{k-1},\eps_ki) \\
&\vdots \\
&(\tfrac12+\tfrac{\sqrt3}2i,\tfrac12-\tfrac{\sqrt3}2i,1,1,i,1,i,\ldots,1,i,\eps_{k-1},\eps_ki)
\end{align*}
in succession.
\end{proof}

Now from Corollary~\ref{cor:connected} we immediately have the following:
\begin{cor}\label{cor:n+2n}
Let $n\in\Nats$, $n\ge2$.
Then $\Fc^\Reals_{n+2,n}$ is connected.
\end{cor}

We take Theorem~\ref{thm:Fk2} and Corollary~\ref{cor:n+2n} to be strong indications:
\begin{conj}
Let $k,n\in\Nats$, with $n\ge2$ and $k\ge n+2$.
Then $\Fc^\Reals_{k,n}$ is connected.
\end{conj}

We also conjecture connectedness in the complex case:
\begin{conj}
Let $k,n\in\Nats$, with $n\ge1$ and $k>n$.
Then $\Fc^\Cpx_{k,n}$ is connected.
\end{conj}

\bibliographystyle{plain}

\end{document}